%% file: hparxiv.tex
\documentclass[11pt]{amsart}
\usepackage{amsopn}
\usepackage{amssymb, amscd}

\usepackage{graphics, epsfig}
\newcommand{\nc}{\newcommand}

\nc{\noi}{\noindent}
\nc{\ft}{\footnote}
\nc{\G}{\Gamma}
\nc{\g}{\gamma}
\nc{\Ld}{\Lambda}
\nc{\ld}{\lambda}
\nc{\la}{\langle}
\nc{\al}{\alpha}
\nc{\be}{\beta}
\nc{\te}{\theta}
\nc{\ra}{\rangle}
\nc{\ba}{\backslash}
\nc{\ke}{\hspace{-.3cm}}
\nc{\msl}{\{\!\{}
\nc{\msr}{\}\!\}}
\nc{\m}{\mathcal}
\nc{\mb}{\mathbf}
\nc{\w}{{\mathbf w}}
\nc{\x}{{\mathbf x}}
\nc{\y}{{\mathbf y}}
\nc{\z}{{\mathbf z}}
\nc{\vv}{{\mathbf v}}
\nc{\bb}{{\mathbf b}}
\nc{\M}{{\mathbf M}}
\nc{\p}[2]{(#1,#2)}

\nc{\SO}{{\mathrm SO}}
\nc{\Spe}{{\mathrm Sp}}
\nc{\Sl}{{\mathrm SL}}
\nc{\SU}{{\mathrm SU}}
\nc{\Or}{{\mathrm O}}
\nc{\U}{{\mathrm U}}
\nc{\Gl}{{\mathrm GL}}
\nc{\Se}{{\mathrm S}}
\nc{\Cl}{{\mathrm Cl}}
\nc{\Spin}{{\mathrm Spin}}
\nc{\Pin}{{\mathrm Pin}}
\nc{\tr}{{\mathrm Tr}}

\nc{\R}{{\mathbb R}}
\nc{\HH}{{\mathbb H}}
\nc{\C}{{\mathbb C}}
\nc{\Z}{{\mathbb Z}}
\nc{\F}{{\mathbb F}}
\nc{\N}{{\mathbb N}}
\nc{\Q}{{\mathbb Q}}
\nc{\PP}{{\mathbb P}}

\nc{\rank}{\operatorname{rank}}

\newtheorem*{theorem}{Theorem}
\newtheorem{lemma}{Lemma}

\title{Hearing the Platycosms}

\author[J.\ P.\ Rossetti \and J.\ H.\ Conway]{J. P. Rossetti$\,^\dag$ \and J. H. Conway$\,^\ddag$}

\address{FaMAF(CIEM), Universidad Nacional de C\'ordoba, 5000-C\'ordoba, Argentina.}
\email{rossetti@mate.uncor.edu}
\address{Department of Mathematics, Princeton University, Princeton, NJ 08544, USA.}
\email{conway@math.princeton.edu}

\thanks{$\dag$ Supported by a Guggenheim fellowship.} 
\thanks{$\ddag$ NSF grant DMS-0072839}

\begin{document}

\maketitle

\noindent{\small \textsc{Note by J.\ H.\ Conway}:
The main result of the following paper is a proof of Rossetti's theorem:
{\it There is, up to scale, a unique isospectral pair of compact platycosms.}

\noindent The proof follows Rossetti's unpublished one (2001) almost exactly step by step --- my main
contribution has been to shorten the exposition, chiefly by introducting conorms to make the argument more transparent,
and by quoting Schiemann's theorem, which Rossetti's proof managed to avoid, except for the torocosm.}

\section{Introduction}

One of John Milnor's great discoveries was that the flat tori based on the Witt lattices
$E_8\oplus E_8$ and $D_{16}^+$ are isospectral. Over the years, the dimension in which isospectral
flat tori were known has been gradually reduced from 16 (Milnor \cite{Mi})
through 12 (Kneser \cite{Kne}), 8 (Kitaoka \cite{Ki}), 6 and 5 (Conway and Sloane \cite{CSi})
to 4 (Schiemann \cite{Schi4}, simplified in \cite{CSi}).  It is now known that that dimension
cannot be further reduced, in view of Schiemann's theorem \cite{Schi3}:

\vskip.15cm

{\it There is no non-trivial isospectrality between flat tori of dimension 3.}

\vskip.15cm

The next simplest candidates are the closed flat manifolds.
They have been discussed in \cite{DM}, \cite{MRp}, \cite{MRl}, where,
among other results, there are some 4-dimensional examples.
On the other hand, it was shown long ago (cf.\ \cite{BGM}, p.\ 153) that there is no non-trivial
isospectrality between flat manifolds in dimension~2.
In our previous paper \cite{CR}, we introduced the term {\it platycosm} (``flat universe'') for
a flat 3-manifold without boundary, since this is the simplest kind of `universe' we might imagine
ourselves to inhabit. In this paper, we close the above gap by proving
\begin{theorem}
There is, up to scale, a unique non-trivial isospectral pair of compact platycosms.
\end{theorem}

This pair was described in \cite{DR}.
The platycosms are described in detail in \cite{CR}, but we have tried to make this paper complete in
itself by repeating some of the preliminary material.

We thank Peter Doyle for discussions on the subject of this paper.
Rossetti thanks Princeton University for its hospitality.

\section{Conorms of Lattices}\label{conorms}

\subsection*{Introduction}
The conorms of a lattice ${\m L}$ are certain numbers that are determined by ${\m L}$ (up to equivalence)
and return the compliment by determining ${\m L}$, at least in low dimensions.
More precisely ${\m L}$ has a {\it conorm function} defined on {\it conorm space}, which is a
finite set that has the structure of a projective $(n-1)$-dimensional space over the field of order two.
At least for $n\le 4$, two $n$-dimensional lattices ${\m L}$ and ${\m L}'$ are isometric if and only if there is an
isomorphism between their conorm spaces that takes one conorm function to the other.

The theory in low dimensions is greatly simplified by the observation that for $n\le 3$ every lattice has
an {\it obtuse superbase}, which implies in particular that the conorms are $\ge 0$.
A {\it superbase} for an $n$-dimensional lattice ${\m L}$ is an $(n+1)$-tuple $\{\mb{v_0},\mb{v_1},\dots,\mb{v_n}\}$ of
vectors that generate ${\m L}$ and sum to zero. It is {\it obtuse} if all inner products $\p {\mb{v_i}}{\mb{v_j}}$ of distinct
vectors are non-positive
(and {\it strictly} obtuse if they are strictly negative).  For a lattice with an obtuse superbase, it can be shown that the
conorms are the negatives of the inner products of pairs of distinct superbase vectors, supplemented by zeros.
(If $n\ge 4$ other things can happen; a lattice may not have an obtuse superbase, and some conorms may be strictly negative.)

{\it For $n=0$, conorm space is empty}, so there are no conorms.

{\it For $n=1$, conorm space is a single point}, and there is one strictly positive conorm $A$.
We represent this by the picture $\,\,\stackrel{A}{\bullet}\,$; it means that the lattice has an obtuse suberbase
$\{\vv,-\vv\}$ with Gram-matrix {\small $\left[\begin{array}{cc} A & -A \\ -A & A \end{array}\right]$};
equivalently a base $\{\vv\}$ with Gram-matrix $\left[ A \right]$.

{\it For $n=2$, conorm space is a 3-point projective line}, and the 3 conorms $A,\,B,\,C$
are non-negative.
We represent this by the picture
$\,\,\stackrel{A}{\bullet}\!\!\!\!-\!\!\!-\!\!\!-\!\!\!\!\stackrel{B}{\bullet}\!\!\!\!-\!\!\!-\!\!\!-\!\!\!\!
\stackrel{C}{\bullet}$~;
it means that the lattice has an obtuse superbase $\{\mb{v_0},\mb{v_1},\mb{v_2}\}$ with matrix 
{\small $\left[\begin{array}{ccc} B+C & -C & -B \\ -C & C+A & -A \\ -B & -A & A+B \end{array}\right]$},
or equivalently it has a base $\{\mb{v_0},\mb{v_1}\}$ with matrix
{\small $\left[\begin{array}{cc} B+C & -C  \\ -C & C+A \end{array}\right]$}.

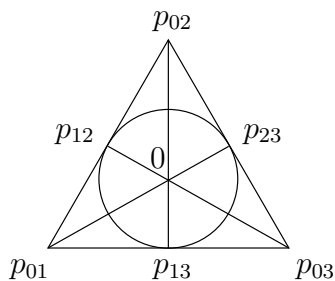
\begin{figure}[h!]
\centerline{\input{fig1.pstex_t}}
\caption{A Fano plane with $p_{01}, p_{02}, p_{03}, p_{12}, p_{13}, p_{23}, 0$ in the usual arrangement}
\label{fig1}
\end{figure}

{\it For $n=3$, conorm space is the 7-point ``Fano plane''},
represented by Figure~\ref{fig1}.
The superbase has Gram-matrix
$$\left[\begin{array}{cccc} p_{0|123} & -p_{12} & -p_{13} & -p_{01} \\
-p_{12} & p_{1|023} & -p_{23} & -p_{02} \\
-p_{13} & -p_{23} & p_{2|013} & -p_{03} \\
-p_{01} & -p_{02} & -p_{03} & p_{3|012}
\end{array}\right],$$
where $p_{ij}=p_{ji}$ are the conorms, whose minimum is $0$, and $p_{i|jkl}:=p_{ij}+p_{ik}+p_{il}$.

The sum of the conorms of all points not in a subspace of codimension 1 is the norm of a (lax) ``Voronoi vector''
\cite{CS7}, that is to say, a vector $\vv$ whose norm $N(\vv)=\p\vv\vv$ is minimal in its coset $(\mathrm{mod}\,\,2{\m L})$.
This is useful because it entails that the conorm function may not be supported on a proper subspace of conorm space,
and so not too many conorms may be zero. In \cite{Co} it is shown, although we shall not need this,
that when $n\le 3$ the location and number $z$ of the zero conorms controls the topology of the Voronoi cell $V$, as follows:

   For $n=1$, $z=0$ and $V$ is an interval.

   For $n=2$, either $z=0$ and $V$ is a hexagon
              or $z=1$ and $V$ is a rectangle.

   For $n=3$, $V$ is topologically:
\begin{itemize}
\item[] a truncated octahedron, if $z=1$;
\item[] a `hexarhombic dodecahedron', if $z=2$;
\item[] a rhombic dodecahedron or hexagonal prism, if $z=3$ \\
\indent according as the three 0s are or are not in line;
\item[]  a cuboid, or `rectangular box', if $z=4$.
\end{itemize}

\subsection*{Putative conorms and the reduction algorithm}

A 3-dimensional lattice has many systems of `putative conorms' in addition to its unique system of  actual
conorms, for instance the numbers obtained by arranging $0$ and the negatives $p_{ij}:=-\p{\mb{v_i}}{\mb{v_j}}$ of the
inner products of distinct members of {\sl any} superbase on a Fano plane in the manner of Figure~\ref{fig1}.
If the putative conorms are all non-negative, they will be the actual conorms.
Otherwise, the following algorithm quoted from \cite{Co} and \cite{CS7} will produce the latter.

Select a `working line' that contains both a $0$ conorm and a negative one, say $-\epsilon$.
Then we transform to an improved system of putative conorms by adding $\epsilon$ to the 3 conorms
on the working line, and subtracting $\epsilon$ from the 4 conorms off this line.
If the improved system still has a negative conorm, we can define a new working line and repeat
the procedure. A finite number of repetitions will suffice to produce the actual conorms.

As an example, Figure \ref{fig0} finds the conorms for the lattice whose Gram-matrix
with respect to a suitable base $\mb{v_1},\mb{v_2},\mb{v_3}$ is
$\left[\begin{array}{ccc} 2&1&1 \\ 1&3&1 \\ 1&1&4 \end{array}\right].$
The Gram-matrix for the superbase $\mb{v_1},\mb{v_2},\mb{v_3},-\mb{v_1}-\mb{v_2}-\mb{v_3}$ is
$\left[\begin{array}{cccc} 2&1&1&-4 \\ 1&3&1&-5 \\ 1&1&4&-6 \\ -4&-5&-6&15 \end{array}\right]$
(found by making each row and column sum to $0$),
which leads to the putative conorms of Figure \ref{fig0}(a).
\begin{figure}[h!]
\input{fig0.pstex_t}
\caption{}
\label{fig0}
\end{figure} 
Transforming this using the working line indicated on it, we obtain Figure \ref{fig0}(b), and we proceed from
this in a similar way to Figures \ref{fig0}(c), \ref{fig0}(d), \ref{fig0}(e),
the last of which gives the actual conorms.
The reader might like to verify that another choice of working line yields the same conorm function.

For a 2-dimensional lattice with putative conorms $A,\,B,\,C$ there is a similar algorithm. We select a
`working point' at which there is a negative conorm, say $-\epsilon$, and transform by adding $2\epsilon$
to this conorm and subtracting $2\epsilon$ from the other two conorms.
Thus
$\,\,\stackrel{-3}{\bullet}\!\!\!\!\!-\!\!\!-\!\!\!-\!\!\!\stackrel{5}{\bullet}\!\!\!-\!\!\!-\!\!\!-\!\!\!
\stackrel{10}{\bullet}\,\,$ transforms through
$\,\,\stackrel{3}{\bullet}\!\!\!-\!\!\!-\!\!\!-\!\!\!\!\!\stackrel{-1}{\bullet}\!\!\!\!\!-\!\!\!-\!\!\!-\!\!\!
\stackrel{4}{\bullet}\,\,$ to
$\,\,\stackrel{1}{\bullet}\!\!\!-\!\!\!-\!\!\!-\!\!\!\stackrel{1}{\bullet}\!\!\!-\!\!\!-\!\!\!-\!\!\!
\stackrel{2}{\bullet}\,$.

Conorms are used in the proof of our theorem in the following way. First, the parameters we use to
specify a particular platycosm $P$ are the conorms of a particular lattice ${\m N}$ associated
to $P$ --- its {\it Naming lattice}. Then, from a supposed isospectrality between $P$ and $P'$ we deduce several
isometries ${\m L}_1\cong {\m L}_1'$, ${\m L}_2\cong {\m L}_2'$,\dots between pairs of lattices associated to $P$ and $P'$.
By equating the corresponding pairs of conorm diagrams we obtain equations in the parameters which if they
do not entail a contradiction show either that $P$ and $P'$ are identical or that they are the unique isospectral tetracosm
and didicosm.

\section{Platycosms and their Universal Covers}\label{platy}

\begin{figure}[!htb]
\centerline{\mbox{\includegraphics*[scale=.4]{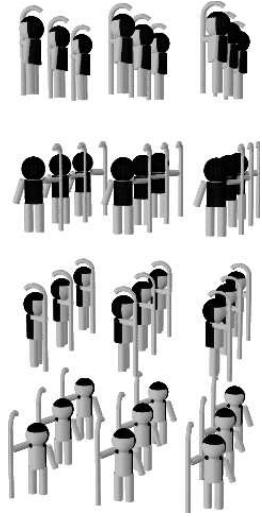}}}
\caption{A tetracosm $c4$.}
\label{fig2}
\end{figure}

\begin{figure}[!htb]
\centerline{\mbox{\includegraphics*[scale=.4]{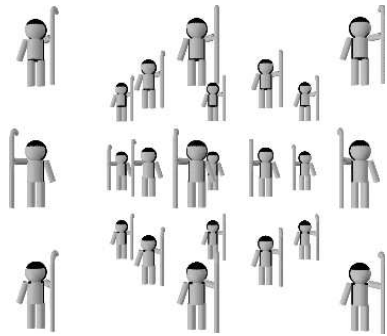}}}
\caption{A second amphicosm $-a1$.}
\label{fig3}
\end{figure}

A locally Euclidean closed 3-manifold is termed
a {\it platycosm} (i.e.\ ``flat universe'').
When you hold something in one hand in a small enough platycosm, you appear to be
surrounded by images of yourself which will either all hold things with the same hand
(if the manifold is orientable, Figure \ref{fig2}) or half with their
left hands and half with their right (if not; Figure \ref{fig3}).

We therefore call a platycosm \,{\it chiral}\, (``handed'') or \,{\it amphichiral}\,
(``either handed'') according as it is or is not orientable.

The images one sees of oneself lie in the manifold's universal cover, which is
of course a Euclidean 3-space $\R^3$.  The symmetry operators that relate them
form a crystallographic space-group $\G$, and geometrically the manifold is the
quotient space $\R^3/\G$.  For this to be a manifold only the identity element
of $\G$ may have a fixed point, a condition that is satisfied by just 10 out of
the 219 three-dimensional space groups.\ft{The number of space groups is often given as 230,
which is inappropiate from the orbifold point of view. The discrepancy arises from the
fact that 11 of the 219 groups arise in two distinct enantiomorphic forms.}

In the notation of \cite{CR}, the 10 (compact) platycosms are:

\

\noindent
\begin{tabular}{r@{\hspace{.3cm}}ccccc @{\hspace{.2cm}} c@{\hspace{.3cm}}cc @{\hspace{.2cm}} c@{\hspace{.2cm}}c}
symbol:      & $c1$  & $c2$  & $c3$  & $ c4$  & $c6$  & $ c22$          & $+a1$  & $-a1$ & $+a2$ & $-a2$ \\
\!\!\!point group: & $C_1$ & $C_2$ & $C_3$ & $ C_4$ & $C_6$ & $C_2\!\times\!C_2$ & $C_2$  & $C_2$ & $C_2\!\times\!C_2$ & $C_2\!\times\!C_2$\\[-.3cm]
\raisebox{-.6cm}[0pt][0pt]{name:}
& \multicolumn{5}{c}{$\underbrace{\phantom{moreC1 C2 C3 C4 C6}}_{\text{\normalsize helicosms}}$} &
$\underbrace{\phantom{C2C2}}_{\text{\normalsize didicosm}}$ &
\multicolumn{2}{c}{\!\!\!$\underbrace{\phantom{C2 C2mm}}_{\text{\normalsize amphicosms}}$} &
\multicolumn{2}{c}{$\underbrace{\phantom{C2xC2 C2xC2}}_{\text{\normalsize amphidicosms}}$}
\\[-.3cm]
\raisebox{-.6cm}[0pt][0pt]{orientable?} & \multicolumn{6}{c}{$\underbrace{\phantom{moreC1 C2 C3 C4 C6 C2xC2mm}}_{\text{\normalsize chiral}}$} &
\multicolumn{4}{c}{\!\!$\underbrace{\phantom{moremmC2 C2 C2xC2 C2xC2}}_{\text{\normalsize amphichiral}}$}
\end{tabular}

\bigskip

Strictly speaking 10 is the number of {\it isotopy classes} of platycosms. The parameters that
determine the metric for a platycosm in any of these classes can be continuously varied so
as to convert it to any other platycosm in the same class. We now discuss these parameters
for a particular platycosm $P$.

Some elements of the space group of $P$ give rise to certain vectors $\vv$; namely $\vv$ may
be obtained from:

a {\it translation} that takes any $\x$ to $\x+\vv$,

a {\it glide reflection} obtained by composing the above translation with some reflection
fixing $\vv$,

a {\it screw motion}, obtained by composing it with some rotation fixing $\vv$.\ft{if this rotation 
has order $N$, we say that the screw motion has {\it period} $N$.}

The {\it Naming lattice} ${\m N}$ of a platycosm is the lattice generated by the vectors so obtained
from all the translations, glide reflections and screw motions of its space group.

Our parameters for $P$ are the conorms of ${\m N}$, arranged in a way described in the appropiate sections
of this article, and more minutely in \cite{CR}.

The translations in $\G$ are usually identified with the corresponding vectors, and form its {\it Translation lattice}
${\m T}$, to which is associated the {\it torus} $T=\R^n/{\m T}$. The group $G:=\G/{\m T}$ is the {\it point group},
which is usually identified with the finite group of orthogonal matrices ${\mathbf M}$ for which some element $\g$ of $\G$ takes
$\x$ to ${\mathbf M}\x+{\mathbf b}$.

\section{What you can hear}\label{hear}

Analytically, the Laplace spectrum $\msl \mu\msr$\ft{We use double braces $\msl\dots\msr$ for multisets --- that
is to say systems of numbers with multiplicities.}
of a manifold is usually coded by the trace
$\sum e^{-\mu t}$ of its heat kernel, and the length spectrum $\msl |\vv|\msr$ of a lattice ${\m T}$ by
its $\theta$-function $\theta_{\m T}(q)=\sum q^{N(\vv)}$. The Appendix shows that these are related by
\begin{equation}\label{trace}
\sum e^{-\mu t} = \frac 1{|G|} \sum_{\g\in G} \frac {\mathrm{vol}({\m T}_{\g})} {(4 \pi t)^{n_{\g}/2}}
\, \theta_{{\m T}_{\g}}(e^{\frac{-1}{4t}}),
\end{equation}
where if $\g$ is represented by $\x\mapsto \M\x+\bb$ we have:
\begin{itemize}
\item[] $n_{\g}$ is the dimension of the fixed space $V_{\M}$ of $\M$,
\item[] ${\m T}_{\g}$ is the orthogonal projection onto $V_{\M}$ of the image of ${\m T}$ under $\g$,
\item[] and the {\it volume} of a lattice is the volume of the associated torus.
\end{itemize}

This is a weighted sum of terms of the form $t^{-a}e^{-b/t}$;
if we put it into the form $\sum c_{a,b} t^{-a}e^{-b/t}$ by collecting terms with the same $a,b$, then it determines
the coefficients $c_{a,b}$. In particular, it determines the subsum
$\sum_{a=d/2}c_{a,b} t^{-a}e^{-b/t}$, which we call its {\it $d$-portion}.
The $n$-portion is $\frac 1{|G|} \frac{{\mathrm{vol}}({\m T})} {(4 \pi t)^{n/2}}
\, \theta_{{\m T}}$ and so since the constant term in $\theta_{\m T}$ is 1, we can hear
$\frac{\mathrm{vol}({\m T})} {|G|}$. But it is also the trace of the heat kernel for $T$
divided by $|G|$, which equals $\frac 1{|G|} \sum_{\mu\in{\rm spec}(T)} e^{-\mu t}$, whose constant
term is $\frac 1{|G|}$, showing that we can also hear $|G|$ and therefore $\mathrm{vol}({\m T})$.
We can now multiply the $d$-portion by $|G|(4\pi t)^{d/2}$ to obtain what we call the {\it $d$-sum} of $M$,
namely, the sum of ${\mathrm{vol}({\m T}_{\g}})\,\theta_{{\m T}_{\g}}$ over those $\g\in G$ for which $n_\g=d$.
Conversely $|G|$ and the $d$-sums determine all the $d$-portions, and so the spectrum.

\smallskip

We summarize:
what we can hear is precisely the order of the point group $G$ together with the $d$-sum for each dimension $d\le n$,
where this is the sum of all the products
$$\mathrm{vol}({\m T}_\g) \theta_{{\m T}_\g}(t)$$
for which $n_\g=d$ is the dimension of the shifted lattice ${\m T}_{\g}$, defined as follows.

Given $\M\in G$ select a vector $\bb$ so that $\g : \x \to \M\x+\bb$ is in the space group $\G$;
then the shifted lattice ${\m T}_\g$ is the projection of $\g({\m T})$ onto the fixed space of $\M$.
This is independent of the choice of $\bb$, since different $\bb$'s will differ by a vector in ${\m T}$,
so we can also call it ${\m T}_{\mb{M}}$.

In particular for $d=n$ the $d$-sum is $\mathrm{vol}({\m T}) \theta_{{\m T}}(q)$, showing again that we can hear
$\mathrm{vol}({\m T})$, since the constant term in $\theta_{\m T}(q)$ is~1.
This proves Sunada's result \cite{Su} that we can hear $\theta_{\m T}(q)$,
which by Schiemann's theorem determines the shape of ${\m T}$ in dimension 3.
We can also hear whether $M$ is orientable, since this happens just if the
$d$-sums vanish for all $d\not\equiv n \mod 2$ (cf. \cite{MRl}).

We remark that the argument shows that we can hear the set of lengths of closed geodesics in $M$,
since these are the lengths of vectors in the various lattices ${\m T}_\g$.  However, we
cannot hear their multiplicities, since the platycosms of our isospectral pair have different numbers
of geodesics of certain lengths (see Section \ref{Sddt}).

\section{The chiral platycosms}

These are the compact orientable flat 3-manifolds without boundary.
They are the {\it helicosms} $c1_{A\,B\,C}^{D\,E\,F}$,
$c2_{A\,B\,C}^{D}$, $c3_{A\,A\,A}^{D}$,
$c4_{A\,A}^{D}$, $c6_{A\,A\,A}^{D}$ and the {\it didicosm}
$c22^{A\,B\,C}$.   The helicosms for $N = 2,3,4,6$ are called the
{\it dicosm, tricosm, tetracosm, hexacosm} respectively.

The {\it torocosm} $c1_{A\,B\,C}^{D\,E\,F}$ is just the 3-dimensional torus based on
the lattice with conorms
\begin{figure}[h!]
\centerline{\input{fig6.pstex_t}}
\label{fig6}
\end{figure}

The space group of any other helicosm $cN$ is generated by the translations of a 2-dimensional
lattice ${\m L}$ that has a rotational symmetry of order $N$, together with a period $N$
screw motion through a vector $\vv$ perpendicular to ${\m L}$.
Its translation lattice ${\m T}$ is generated by ${\m L}$ together with $N\vv$.
The lower parameters (typically $A\,B\,C$) are the (non-zero) conorms of ${\m L}$, while
the upper parameter $D$ is the norm of $\vv$.

\begin{figure}[h!]
\centerline{\mbox{\includegraphics*[scale=.65]{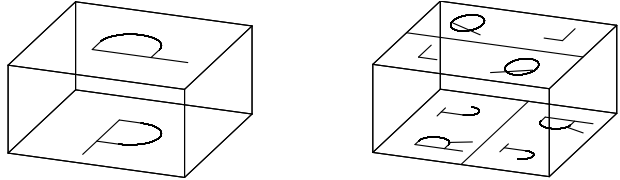}}}
\caption{Tetra (left) and Didi (right) are obtained by identifying the sides by translations, and top and bottom
as shown.}
\label{fig7}
\end{figure}

The space group of the {\it didicosm} $c22^{A\,B\,C}$ is generated by three period 2 screw
motions whose associated vectors are perpendicular, of norms $A,B,C$. Its translation lattice ${\m T}$
is generated by the doubles of these vectors.
For a more complete description, see \cite{CR}, or Figures \ref{fig7} and \ref{fig8} in Section \ref{Sddt}.

\subsection*{Chiral platycosms of the same type}
For manifolds both of type $cN$, the result is due to Schiemann \cite{Schi3} when $N=1$.
If $N>1$, the generic form
is $cN_{A\,B\,C}^D$. Now for this we can hear the multiset of numbers
$\msl A,B,C,N^2D\msr$ since these are the conorms of ${\m T}$ after
some $0$s have been deleted.  But we can also hear $N$ (the order of $G$) and $D$
(the minimal norm in the 1-sum), whence also $N^2D$ and therefore
$\msl A,B,C\msr$ by deletion.  But now $N,D$ and $\msl A,B,C\msr$
determine the manifold.

For manifolds both of type $c22$ the argument is even easier. We can
hear $\msl A,B,C\msr$ since we can hear $\msl4A,4B,4C\msr$,
the conorms of ${\m T}$, and these determine $c22^{A\,B\,C}$ since this is invariant under all 
permutations of $A,B,C$.

\subsection*{Chiral platycosms of different types}
Since the point groups of the six chiral platycosms have the respective orders 1, 2, 3, 4, 6, 4,
an isospectral couple that are not of the same type
must consist of a didicosm $c22^{A\,B\,C}$ say, and a tetracosm $c4_{E\,E}^{D}$.
Moreover, we must have $\msl 4A,4B,4C\msr=\msl 16D,E,E\msr$, since these
are the non-zero conorms of the two translation lattices,
so we can suppose the isospectrality is $c22^{A\,\,4D\,A}\sim
c4_{4A\,\,4A}^{D}$.
Now let $A=a^2$, $D=d^2$.
Then the 1-sums are
{\setlength\arraycolsep{1pt}
\begin{displaymath}
\begin{array}{lrcccl}
& 2a\,\theta\msl A \text{ odd}^2\msr_A &+& 4d\,\theta\msl 4D\text{ odd}^2 \msr_{4D} &+&
2a\,\theta\msl A \text{ odd}^2\msr_A \qquad \nonumber \\
\text{and}\quad & \nonumber \\
& 4d\,\theta\msl D {(4n+1)}^2\msr_D &+& 4d\,\theta\msl D {(4n+2)}^2  \msr_{4D} &+&
4d\,\theta\msl D {(4n+3)}^2\msr_D \nonumber
\end{array}
\end{displaymath}}
where ``$\theta\msl\dots\msr_K$'' denotes the $\theta$-function of a lattice whose norms
are the multiset $\msl\dots\msr$, which has minimal non-zero member $K$.
Then from the equality of minimal norms here we deduce $A=D$, completely
determining the putative isospectrality:
$$c22^{D\,4D\,D}\sim c4_{4D\,4D}^{D}.$$
This is indeed an isospectrality since the middle terms in
{\setlength\arraycolsep{1pt}
\begin{displaymath}
\begin{array}{lrcccl}
 & 2d\,\theta\msl D \text{ odd}^2\msr_D &+& 4d\,\theta\msl 4D\text{ odd}^2 \msr_{4D} &+&
2d\,\theta\msl D \text{ odd}^2\msr_D \\
\text{and} & \qquad \\
 & 4d\,\theta\msl D {(4n+1)}^2\msr_D &+& 4d\,\theta \msl D {(4n+2)}^2 \msr_{4D} &+&
4d\,\theta\msl D {(4n+3)}^2\msr_D \nonumber
\end{array}
\end{displaymath}}
are equal, while the outer ones combine correctly using the fact that each
odd number is of just one form $4n+1$ or $4n+3$.

\section{The DDT-example}\label{Sddt}

Since it involves a Didicosm and a Tetracosm, Doyle has called this the DDT-example \cite{DR}.
We have proved the {\it DDT-theorem} announced there:
{\it `Didi', the didicosm $c22\framebox{$^{1\,1\,2}$}$,
is isospectral to `Tetra', the tetracosm $c4\framebox{$^1_{2\,2}$}$.}\ft{The convention
that
\framebox{$_{a\,b\,c}^{d\,e\,f}$}
means
$_{a^2\,b^2\,c^2}^{d^2\,e^2\,f^2}$
was introduced in \cite{CR}. The numbers inside these two squares are the lengths of
the three orthogonal translations that generate the naming lattices of Tetra and Didi.}

\begin{figure}[h!]
\centerline{\mbox{\includegraphics*[scale=.65]{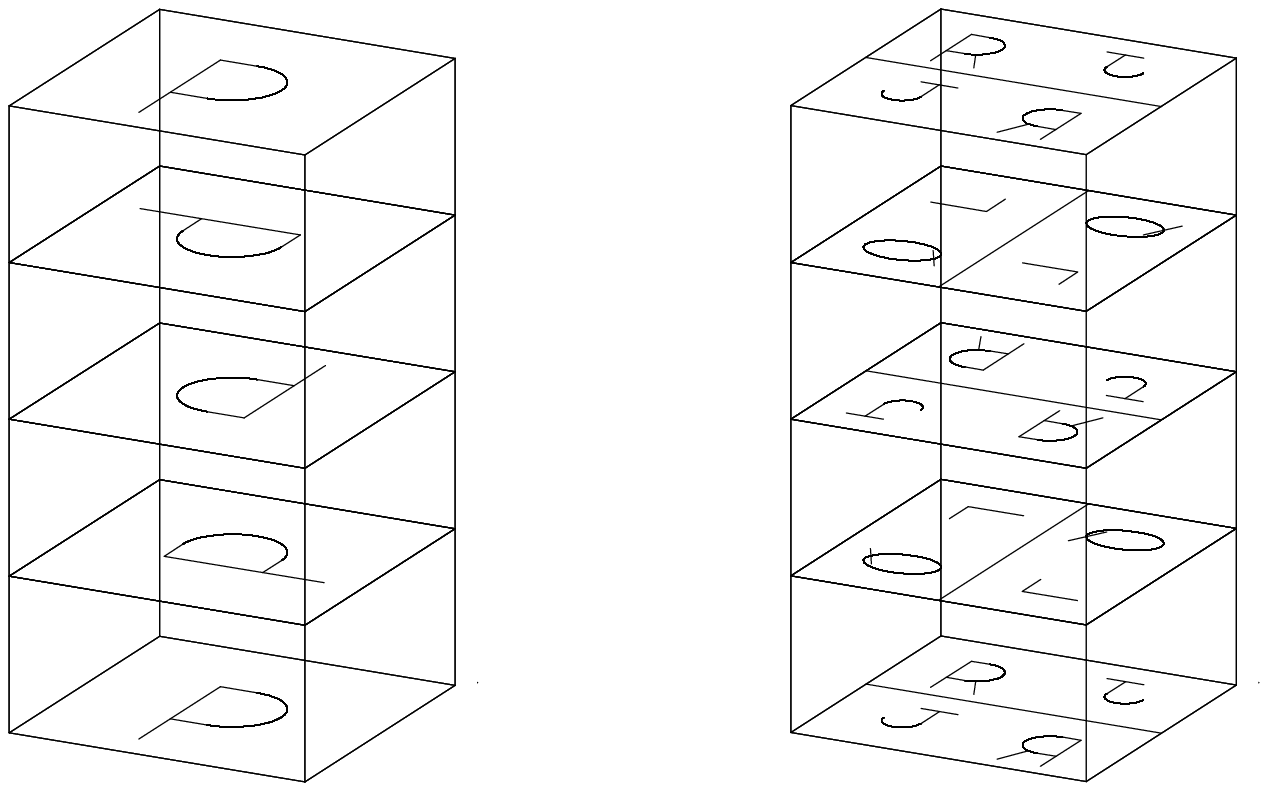}}}
\caption{}
\label{fig8}
\end{figure}

We have doubled the scale used in \cite{DR} so as to make the minimal geodesics have length 1. As a mnemonic,
the most natural fundamental region for Tetra is a box of volume 4, while for Didi it is really the union of
2 boxes of volume 2.
The $2\times 2\times 4$ boxes of Figure \ref{fig8}  are fundamental regions for the associated (isometric) translation lattices
${\m T}$ and ${\m T'}$.
It is clear that the shortest geodesics have length 1, but in Didi they are four
`half-turn' ones, in two different directions, while Tetra has only two parallel `quarter-turn' ones. This
shows

{\it Although we can hear the set of geodesic lengths, we cannot hear the multiplicities to which these lengths appear.}

\section{Amphicosms: the balanced case}

\begin{figure}[h]
\centerline{\input{fig11.pstex_t}}
\caption{The first (or positive) and second (or negative) amphicosms}
\label{fig11}
\end{figure}
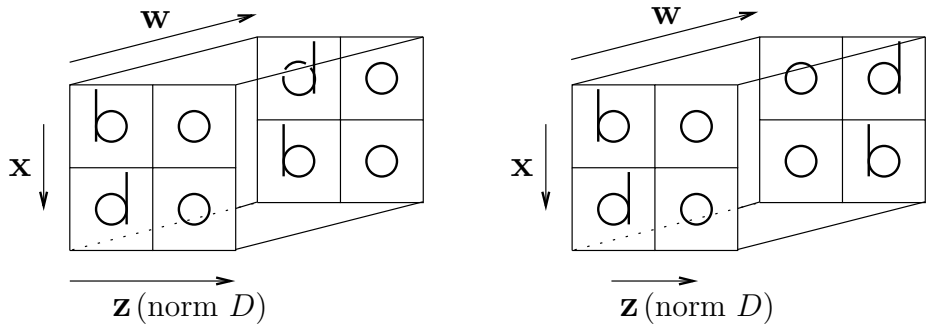

The space group of either of the two {\it amphicosms} $\pm a1^D_{A:B\,C}$ is generated by 
two glide reflections whose associated vectors are $\w$ and $\x$, not necessarily perpendicular,
together with a translation perpendicular to them whose vector we call $\z$ for $+a1$ and $2\z$ for $-a1$.
The two glide reflections are those that take the leading ``{\bf b}'' to the two ``{\bf d}''s in the corresponding
pictures (Figure \ref{fig11}).
The reflecting planes are identical for the first amphicosms, but differ by $\frac 12 \z$ for
the second ones.
$D$ is the norm of $\z$ and $A:B\,C$ the conorms of $\langle \w,\x\rangle$, the conorm $A$ separated by the colon
being $-\p\w\x$.

The point group of the amphicosms has order 2, being generated by a single reflexion $R$.
The equations of $d$-sums therefore take the form
$$\mathrm{vol}({\m T})\,\theta_{\m T}=\mathrm{vol}({\m T}')\,\theta_{{\m T'}} \quad\hbox{ and }
\quad \mathrm{vol}({\m T}_R)\,\theta_{{\m T}_R}=\mathrm{vol}({\m T}_{R'})\,\theta_{{\m T}_{R'}}.$$

The {\it balanced case} handled in this section is when the two shifted lattices ${\m T}_R$ and ${\m T}_{R'}$ have the
same volume, so we have $\theta_{{\m T}_R}=\theta_{{\m T}_{R'}}$.

Recall that the shifted lattice ${\m T}_R$ is the projection onto the fixed space of $R$ of the image of ${\m T}$
under a suitable element $\g$ of $\G$.
Now for the first amphicosm $+a1$, ${\m T}$ is spanned by $2\w,2\x,\w+\x,\z$, while a suitable $\g$ has $\w$ for
its translation part, so that  ${\m T}_R=\langle 2\w,2\x,\w+\x\rangle+\w=\{m\w+n\x:m\not\equiv n\mod 2\}$.
A similar calculation shows that in fact ${\m T}_R$ has the same form also for the second amphicosm $-a1$.
So in either case ${\m T}_R$ is the
difference of a `double' lattice
${\m D}:=\{m\w+n\x: \hbox{ all } m,n\in\Z\}$ and a `half' lattice
${\m H}:=\{m\w+n\x:m\equiv n\mod 2\}$
and its $\theta$-function is $\theta_{\m D}-\theta_{\m H}$.

\begin{lemma}\label{lemma1}
If ${\m D},{\m H},{\m D'},{\m H'}$ are four 2-dimensional lattices of which ${\m H}$
has index 2 in ${\m D}$ and ${\m H'}$ has index 2 in ${\m D}'$,
then $\theta_{\m D}-\theta_{\m H}=\theta_{{\m D'}}-\theta_{{\m H'}}$ holds if and only if ${\m D}'$ and ${\m H'}$ are
respectively isometric to ${\m D}$ and ${\m H}$.
\end{lemma}
\begin{proof}
To say that the double lattice ${\m D}=\langle\mb{v_0},\mb{v_1},\mb{v_2}\rangle$ has conorms $A\,B\,C$ entails that
its typical vector has norm
$Am^2+Bn^2+Cp^2$ where $m+n+p=0$, and we can suppose the vectors of ${\m D}\smallsetminus H$ are those with $m$ odd.
For such a vector $\vv$, if $n=0$, $\vv$ is an odd multiple of $\mb{v_1}$, whose norm is $A+C$;
if $p=0$, $\vv$ is an odd multiple of $\mb{v_2}$, whose norm is $A+B$;
while all other norms in ${\m D}\smallsetminus{\m H}$ are at least $A+B+C$.

Now we can hear the minimal norm of ${\m D}\smallsetminus{\m H}$, say $A+B$, and then we can hear the next primitive
norm $A+C$, since this will be the minimal norm after subtracting $\theta\msl (A+B)n^2 : n\, odd \msr_{A+B}$.

Finally, since the determinant\ft{By the {\it determinant} of a lattice we mean the determinant of its Gram-matrix,
equivalently the squared volume of the associated torus.} of ${\m D}$ is $AB+BC+AC=(A+B)(A+C)-A^2$, the equality
of volumes tells us that we can hear $A^2$, therefore $A$ itself since $A\ge 0$, and thence also $\msl B,C\msr$.
\end{proof}

This lemma encodes an isospectrality between two amphicosms $M,M'$ in three lattice isometries
\,${\m T}\cong {\m T}'$, ${\m D}\cong {\m D'}$, ${\m H}\cong {\m H'}$ (since $\mathrm{vol}({\m T})=\mathrm{vol}({\m T}')$), 
which is important for our strategy because such isometries are easily handled using conorms.
We now discuss the conorms of the three lattices for each type of amphicosm.

For the first amphicosm $+a1$ the naming lattice ${\m N}=\langle \w,\x,\z\rangle$ is spanned by
two glide vectors $\w,\x$ together with 
the translation vector $\z$ where the vectors $\w,\x,\y:=-\w-\x$ and $\z$ have Gram-matrix

\begin{displaymath}
\left[
\begin{array}{cccc}
         A+B  &  -A  &  -B  &  0 \\
         -A  & A+C & -C  & 0 \\
         -B  & -C  & B+C & 0 \\
         0  & 0  & 0  & D
\end{array}
\right] \qquad (A,B,C,D \ge 0),
\end{displaymath}
and since the `double' lattice $\m D$ is spanned by $\w$ and $\x$, we have
${\m D}=\Ld_{A\,B\,C}$, the 2-dimensional lattice with conorms $A\,B\,C$.

What are the conorms of \,${\m H}=\langle 2\w,2\x,\y\rangle$?
One superbase for ${\m H}$ is \,\,$2\w$, $\y$, $-2\w-\y$, whose Gram-matrix is
\begin{displaymath}
\left[
\begin{array}{cccc}
   4A+4B  &  -2B  &  -4A-2B \\
   -2B  & B+C & B-C \\
   -4A-2B  & B-C & 4A+B+C
\end{array}
\right].
\end{displaymath}
The conorms of ${\m H}$ will be the negatives $4A+2B,\,2B,\,C-B$ of the inner products here unless
$C<B$, when they will be the numbers $4A+2C,\,2C,\,B-C$ obtained similarly from the alternative superbase
$2\x,\y,-2\x-\y$.

Since it often happens that the conorms of a lattice take one of several different forms in this way, we introduce
some useful terminology. Either triple
$$4A+2B,\,2B,\,C-B \quad {\rm or} \quad 4A+2C,\,B-C,\,2C$$
is a set of putative conorms
that determines the lattice ${\m H}$, and
this pair of triples is
an {\it exhaustive system}, in the sense that at least one of them
will be the actual conorms.  Finally, we introduce the special notation $[A]\,B\,C$ for the actual
conorms\ft{Beware: the triple of numbers called ${[A]\,B\,C}$ probably contains neither $B$ nor $C$!} here, so that
\begin{displaymath}
\Ld_{[A]\,\,B\,\,C} \quad\text{means}\quad\left\{ \begin{array}{ll}
\Ld_{4A+2B\,\,\,2B\,\,\,C-B} & \textrm{ if } C\ge B; \\
\Ld_{4A+2C\,\,\,B-C\,\,\,2C} & \textrm{ if } B\ge C.
\end{array}\right.
\end{displaymath}

The three sublattices of index 2 in $\Ld_{A\,B\,C}$ are $\Ld_{[A]\,\,B\,\,C}$, $\Ld_{A\,\,[B]\,\,C}$ and
$\Ld_{A\,\,B\,\,[C]}$ and their intersection is $\Ld_{4A\,\,4B\,\,4C}$ (see Figure \ref{fig12}).

\begin{figure}[h!]
\centerline{\input{fig12.pstex_t}}
\caption{}
\label{fig12}
\end{figure}

In fact ${\m D}$ and ${\m H}$ take these forms also for the second amphicosm.
The fact that ${\m H}=\Ld_{[A]\,B\,C}$ rather than $\Ld_{A\,[B]\,C}$ or $\Ld_{A\,B\,[C]}$
distinguishes the separated parameter $A$ in our names $+a1^D_{A:B\,C}$ and $-a1^D_{A:B\,C}$
from $B$ and $C$.
On the other hand, $B$ and $C$ can be interchanged without affecting either manifold.

\hspace{-2pt}In this notation the first amphicosm $+a1^D_{A:B\,C}$ whose naming lattice $\langle \w,\x,\y,\z\rangle$
is $\Ld^D_{A\,B\,C}$
has translation lattice ${\m T}=\langle 2\w,2\x,\y,\z\rangle$ of type $\Ld^D_{[A]\,B\,C}$
(meaning either $\Ld^D_{4A+2B\,2B\,C-B}$ or $\Ld^D_{4A+2C\,B-C\,2C}$).
The second amphicosm with the same naming lattice has
${\m T}=\langle 2\w,2\x,2\y,\y+\z\rangle$, for which an exhaustive system of putative conorms is given in Figure \ref{fig13}.

\begin{figure}[h!]
\centerline{\input{fig13.pstex_t}}
\caption{}
\label{fig13}
\end{figure}

Two isospectral amphicosms (of any signs) must have the same multiset $\msl A,B,C\msr$
since each has ${\m D}=\Ld_{A\,B\,C}$, but their separated parameters might be distinct.
If so, we can suppose they are $\pm a1_{A:B\,C}^?$ and $\pm a1_{B:A\,C}^?$ with $A\ne B$,
and so
$$\text{either} \quad \msl 4A+2B,\,2B,\,C-B\msr \quad \text{or} \quad \msl 4A+2C,\,B-C,\,2C \msr $$
\noindent must equal
$$\text{either} \quad \msl 4B+2A,\,2A,\,C-A\msr \quad\text{or} \quad \msl 4B+2C,\,A-C,\,2C\msr,$$
since these are the conorms of the two half lattices.
But equating the sum of conorms for the two left candidates gives
$$4A+3B+C=3A+4B+C \quad \hbox{and so} \quad A=B,$$
while equating the right candidates similarly gives
$$4A+B+3C=A+4B+3C \quad \hbox{and so again} \quad A=B,$$
contradicting our assumption.

Otherwise we must equate a left candidate with a right one, say
$$4A+2B,2B,C-B \quad \hbox{with} \quad 4B+2C,A-C,2C.$$
Then $C-B$ must equal $A-C$, since it is smaller than $4B+2C$ and $2C$,
and from the sum of the remaining pairs of conorms we obtain $4A+4B=4B+4C$, whence $A=C$, whence $C=B$,
whence again $A=B$, contradicting our assumption.

We are left with the case in which the separated parameters are equal.
Then the two manifolds will be the same if their signs are, since $D$ can be obtained
from the determinant.
Otherwise we can take the isospectrality to be between
$+a1_{A:B\,C}^{4D}$ and $-a1_{A:B\,C}^D$ and suppose $B\le C$.

The conorms
\begin{figure}[h!]
\centerline{\input{fig14.pstex_t}}
\caption{Conorms of ${\m T}$ for $+a1_{A:B\,C}^{4D}$}
\label{fig14}
\end{figure}
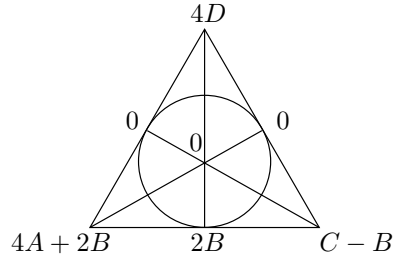
of the translation lattice ${\m T}$ of the first amphicosm (Figure~\ref{fig14}) have the property
\begin{displaymath}
\begin{array}{lc}
(0^3)\,\, : \,\,\, &  \text{ there are three non-collinear zeros, and a fourth zero if $B=0$.}
\end{array}
\end{displaymath}
So this must also be true of the conorms of the translation lattice ${\m T}$ of the second amphicosm, which are one of those
in Figure~\ref{fig15} (simplifying Fig. \ref{fig13}).
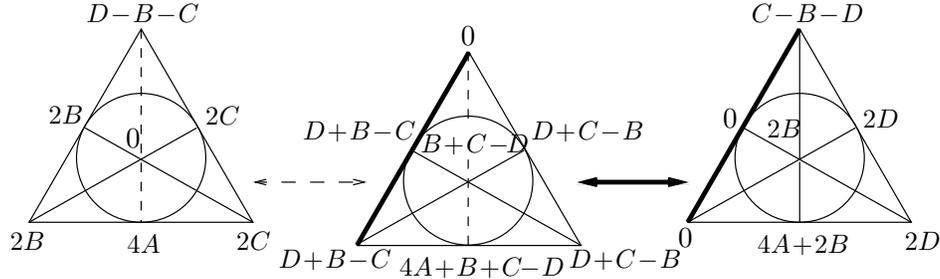
\begin{figure}[h!]
\centerline{\input{fig15.pstex_t}}
\caption{Conorms of ${\m T}$ for $-a1_{A:B\,C}^D$}
\label{fig15}
\end{figure}

Now $D$ cannot vanish since the determinant is $16D(AB+BC+AC)$, and for the same reason $C$ cannot
vanish (since then $B$ would).
If the middle triangle has 3 non-collinear zeros then $D+B-C$ and $B+C-D$ must vanish, since
$4A+B+C-D$ and $D+C-B$ exceed them.
If the left or right one does, we must have $B=0$, requiring a fourth zero, which can only be
the top number.

But all the 3 cases give $B=0,\,D=C$, for which the two translation lattices~${\m T}$ and~${\m T}'$
are in fact distinct (Figure~\ref{fig16}).
\begin{figure}[h!]
\centerline{\input{fig16.pstex_t}}
\caption{}
\label{fig16}
\end{figure}

\section{Amphicosms: the unbalanced case}

If the two volumes in
$$\mathrm{vol}({\m D}\smallsetminus {\m H}) \,\theta_{{\m D}\smallsetminus {\m H}}=
\mathrm{vol}({\m D'}\smallsetminus {\m H'}) \,\theta_{{\m D'}\smallsetminus {\m H'}}$$
are distinct, we may suppose \,$\theta_{{\m D}\smallsetminus {\m H}}=
f \theta_{{\m D'}\smallsetminus {\m H'}}$\, for some \,$f>1$.
But ${\m D}'\smallsetminus {\m H'}$ has at least one pair $\pm\mb{u}$ of minimal vectors,
so ${\m D}\smallsetminus {\m H}$ must have
at least two such pairs, say $\pm \vv$, $\pm \w$.

\begin{lemma} Under these conditions we can assert:
\begin{itemize}
\item[(i)]  $f=2$.
\item[(ii)]  ${\m D}=\langle \vv,\w\rangle$ has type $\Ld_{\alpha\,\beta\,\beta}$.
\item[(iii)] ${\m H}=\langle 2\vv,2\w,\vv+\w\rangle$ has type $\Ld_{[\alpha]\,\beta\,\beta}$.
\item[(iv)] ${\m D}'$ has type $\Ld_{\alpha\,\beta\,[\beta]}$.
\item[(v)]  ${\m H'}$ has type $\Ld_{4\alpha\,4\beta\,4\beta}$.
\end{itemize}
\end{lemma}
\begin{proof}
The double lattice ${\m D}$ certainly contains $\vv$ and $\w$, so is at least the rhombic lattice in Figure~\ref{fig17},
while the shifted lattice ${\m D}\smallsetminus {\m H}$ is at least the rectangular one.
If either were strictly bigger, then every rectangle of the figure would contain a further point in each of
${\m H}$ and ${\m D}\smallsetminus {\m H}$, a contradiction since those in the rectangle around $0$ would be shorter
than $\vv$ and $\w$.

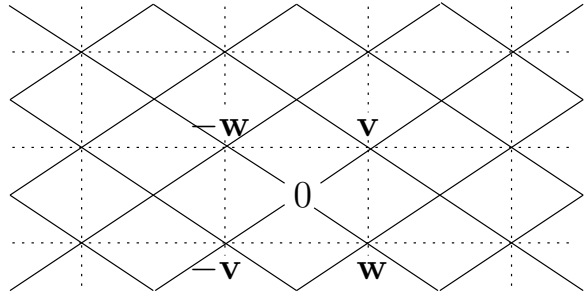
\begin{figure}[h!]
\centerline{\input{fig17.pstex_t}}
\caption{The rhombic and rectangular lattices}
\label{fig17}
\end{figure}

This establishes the structure of ${\m D}$ and ${\m H}$, and the fact that $f=2$.
Now look at the index 2 sublattices of ${\m D}$ and their intersection (Figure~\ref{fig18}$(a)$).

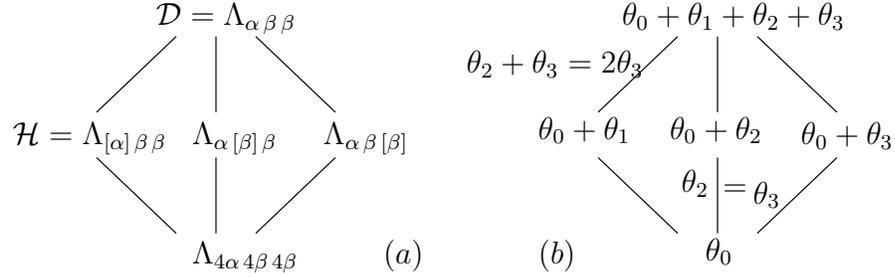
\begin{figure}[h!]
\input{fig18.pstex_t}
\caption{Several lattices and their $\theta$-functions.}
\label{fig18}
\end{figure}

Their $\theta$-functions are as in Figure \ref{fig18}$(b)$, where the $\theta_i$ are those of the 4 cosets of
$\Ld_{4\alpha\,4\beta\,4\beta}$. The isometry  $\Ld_{\alpha\,[\beta]\,\beta}\cong \Ld_{\alpha\,\beta\,[\beta]}$
implies $\theta_2=\theta_3$, showing that  $\Ld_{\alpha\,\beta\,[\beta]}$ and  $\Ld_{4\alpha\,4\beta\,4\beta}$
are possible choices for ${\m D}'$ and ${\m H'}$, since for them
$$\theta_{{\m D'}}-\theta_{{\m H'}}=\theta_3, \quad \theta_{{\m D}}-\theta_{{\m H}}=\theta_2+\theta_3=2\theta_3.$$
Any other pair ${\m D},{\m H}$ must yield the same $\theta_{{\m D}}-\theta_{{\m H}}$, and so, by Lemma 1, be isometric 
to this pair.
\end{proof}

We now know the lower three parameters $A\,B\,C$ for each amphicosm using ${\m D}=\Ld_{A\,B\,C}$, and which of these
is the separated one using ${\m H}=\Ld_{[A]\,B\,C}$.
We can relate their upper parameters $D$ using the determinants of their translation lattices, namely
$$4D(AB+BC+AC) \,\,\, {\rm for} \,\, +a1_{A:B\,C}^{D} \hbox{ and } 16D(AB+BC+AC) \,\,\, {\rm for} \,\,\, -a1_{A:B\,C}^D.$$

In this way we find that the isospectrality must be between one of
$$ +a1_{\alpha:\beta\,\beta}^{16\delta}  \,\,(M_+)  \quad \text{or} \quad  -a1_{\alpha:\beta\,\beta}^{4\delta}
 \,\,(M_-)\qquad\quad\phantom{nadanad}$$
\noindent and one of
$$ +a1_{\alpha-\beta:2\beta\,6\beta}^{4\delta} \,\,(M_+')  \quad \text{or} \quad
 -a1_{\alpha-\beta:2\beta\,6\beta}^{\delta}  \,\,(M_-')  \qquad\quad (\alpha\ge\beta)$$

$$ +a1_{\beta-\alpha:2\alpha\,2\alpha+4\beta}^{4\delta} \,\,(M_+')  \quad \text{or} \quad
 -a1_{\beta-\alpha:2\alpha\,2\alpha+4\beta}^{\delta}  \,\,(M_-')  \qquad\quad (\alpha\le\beta)$$

Putting these values into Figures \ref{fig14} and \ref{fig15}, we obtain the possible conorms
for the corresponding ${\m T}$ and ${\m T}'$, given in Figures~\ref{fig19} and~\ref{fig20} respectively,
\begin{figure}[h!]
\input{fig19.pstex_t}
\caption{${\m T}$ is one of these two lattices.}
\label{fig19}
\end{figure}
where for definiteness we have rescaled to make $\beta=1$ ($\beta$ cannot vanish
since $\alpha\beta+\beta\beta+\beta\alpha$ divides the determinant).
Those for ${\m T}_+$ and ${\m T}_-$ have the respective properties:
\begin{displaymath}
\begin{array}{ll}
(0^4) \,: \,\,& \text{ there are 4 zeros.} \\
(*^4) \,: & \text{ the 4 conorms off some line $l$ are equal and not zero}.
\end{array}
\end{displaymath}

\begin{figure}[h!]
\input{fig20.pstex_t}
\caption{${\m T}'$ is one of these two lattices.}
\label{fig20}
\end{figure}

\begin{figure}[h!]
\input{fig21y22.pstex_t}
\caption{Condition $(0^4)$ yields these.}
\label{fig21y22}
\end{figure}
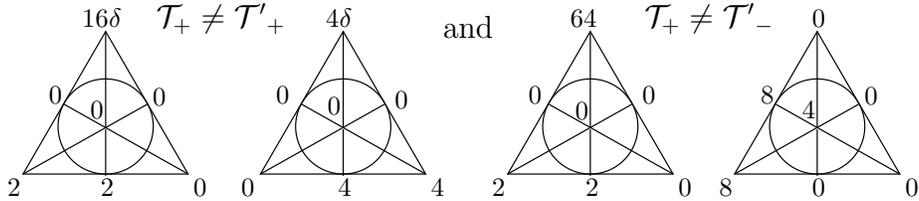

The property $(0^4)$ holds for $M_+'$ only if $\alpha=0$, but the resulting lattices (see Figure \ref{fig21y22} left)
cannot be isometric.

By inspection $(0^4)$ holds for $M_-'$ only if $\alpha=0,\delta=4$, but again the resulting lattices
(${\m T}$ and ${\m T}'$ in Figure \ref{fig21y22} right) are different.

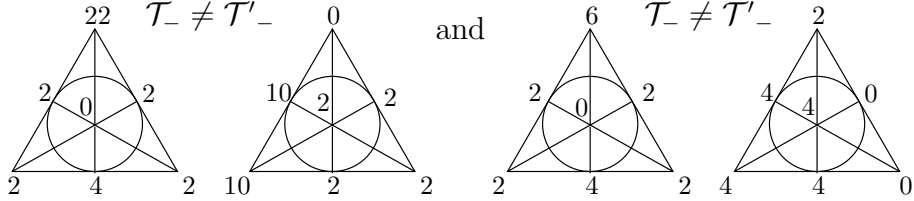
\begin{figure}[h!]
\input{fig23.pstex_t}
\caption{Condition $(*^4)$ yields these.}
\label{fig23}
\end{figure}

Finally property $(*^4)$ cannot hold for $M_+'$ and
holds for $M_-'$ only in the two cases $\alpha=1,\delta=6$ and $\alpha=1, \delta=2$ of ${\m T'}_-$,
but each of the resulting pairs (see Figure \ref{fig23})
also consists of two distinct lattices.

\section{The amphidicosms}

The amphichiral platycosms whose point group has order four are the {\it first} (or {\it positive}) and 
{\it second} (or {\it negative}) {\it amphidicosms}.
Here the argument is simple because the naming lattice ${\m N}=\langle \x,\y,\z\rangle$ is spanned by three
orthogonal vectors.

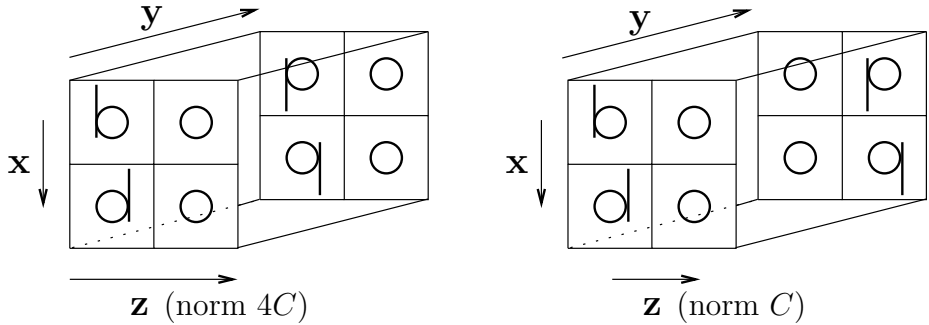
\begin{figure}[h!]
\input{fig24.pstex_t}
\caption{The first amphidicosm $+a2^{4C}_{A:B}$ (left) and the second amphidicosm $-a2^{C}_{A:B}$ (right).}
\label{fig24}
\end{figure}

This case is unusual in that the 2-sum involves two distinct shifted lattices, say ${\m D}_1\smallsetminus {\m H}_1$ and
${\m D}_2\smallsetminus {\m H}_2$.

In view of the anomaly between the meaning of $\z$ in the two cases (see Figure~\ref{fig24}), we shall take
our platycosms in the forms
$$+a2^{4C}_{A:B} \,\,\text{\bf or}\,\, -a2^C_{A:B}$$
so that their tori have the same volume $8\sqrt{ABC}$. Then the 2-sum is given by one of the two formulae
\begin{displaymath}
4\sqrt{AB}\,\theta\msl A .odd^2 + B .even^2 \msr_A +
4\sqrt{BC}\,\theta\msl B .odd^2 + C .(even^2 \,{\bf or}\, odd^2) \msr_{B \,\,{\bf or }\,\, B+C}
\end{displaymath}
respectively, where again the subscripts on the multisets indicate their least elements.

The colons in the notations remind us that this time there is no symmetry --- the three parameters
$A,B,C$ play distinct roles.  However, $B$ is audible since the 1-sum is
\begin{displaymath}
2\sqrt{B}\,\theta\msl B .odd^2\msr_B
\end{displaymath}
and then the multiset $\msl A,B,C \msr$, whence $\msl A,C \msr$, is audible since the 3-sum is
\begin{displaymath}
8\sqrt{ABC}\,\theta\msl 4A \, m^2 + 4B\, n^2 + 4C \, p^2 \msr \qquad m,n,p\in\Z.
\end{displaymath}

The only nontrivial putative isospectrality between amphidicosms with the same $A$ and the same $C$ is between
$+a2^{4C}_{A:B}$ and $-a2^C_{A:B}$, which is easily disproved since the elements that remain after cancelling the first
summands of the two formulae, namely
\begin{displaymath}
4\sqrt{BC}\,\theta\msl B .odd^2 + C .(even^2 \,{\bf or}\, odd^2) \msr_{B \,\,{\bf or }\,\, B+C}\,,
\end{displaymath}
will have distinct minima $B$ and $B+C$.
So we may suppose the isospectrality is between
$+a2^{4C}_{A:B}$ or $-a2^{C}_{A:B}$ and $+a2^{4A}_{C:B}$ or $-a2^{A}_{C:B}$, and, by symmetry, that $A<C$.
But then we must equate one of the original two formulae to
\begin{displaymath}
4\sqrt{BC}\,\theta\msl C .odd^2 + B .even^2 \msr_C +
4\sqrt{AB}\,\theta\msl B .odd^2 + A .(even^2 \,{\bf or}\, odd^2) \msr_{B \,\,{\bf or }\,\, A+B}
\end{displaymath}
the $C$ in which must therefore equal $B$ since it cannot be $A$ or $B+C$.
The second summand of the original formula is now known to be
\begin{displaymath}
4\sqrt{B^2}\,\theta\msl B .odd^2 + B .even^2 \msr_{B},
\end{displaymath}
and coincides with the first summand of the new one.
After cancelling them the minima become
$$A \,\,\text{ and } \,\,\big( B \text{ \bf or } A+B\big),$$
and so $A$ must also equal $B$ since it cannot equal $A+B$. This contradicts our assumption that $A<C$,
and finishes the proof.

\section*{Appendix: A Jacobi-Poisson summation formula for flat manifolds}

We first exhibit some eigenfunctions on the torus $T=\R^n/\m{T}$ defined by
a given lattice $\m{T}$. These are the functions $f_w(x):=e^{2\pi i \p wx}$, for each
$w=(w_1,\dots,w_n)$ in the dual lattice $\m{T}^*:=\{w:\p wv\in\Z\text{ for each } v\in\m{T}\}$.
For since $f_w(x)$ is unaltered if we increase $x$ by some $v\in\m{T}$, it is indeed
defined on $T$, and since differentiating with respect to $x_k$ multiplies $f_w(x)$ by
$2\pi i$ we see that
$$
\left(\frac{\partial^2}{{\partial x_1}^2}+\dots+\frac{\partial^2}{{\partial x_n}^2}\right) f_w(x) =
(2\pi i)^2 \left(w_1^2+\dots+w_n^2\right) f_w(x),$$
showing that $f_w$ is an eigenfunction of $\nabla^2$ with eigenvalue $-4\pi^2 \p ww$.
In fact the Stone-Weierstrass theorem shows the $f_w$ to be a complete system of
eigenfunctions so that the spectrum of $\nabla^2$ is the multiset $-4\pi^2\msl N(w):w\in\m{T}^*\msr$,
where we write $N(w)=\p ww$.

Now each of the flat manifolds $M$ we are concerned with is the quotient of a torus $T$ by
a finite group $G$ (that we can identify with the point group). In these circumstances, every
eigenfunction on $M$ can be regarded as an eigenfunction on $T$, and conversely, an eigenfunction
on $T$ is one on $M$ just if it is invariant under $G$. Let $p$ denote the projection operator
that takes any function with eigenvalue $\mu$ to the mean of its images under $G$ --- then
the multiplicity $m=m_\mu$ of $\mu$ will be the trace of $p$, since $p$ takes the form diag$(1^m,0^{n-m})$
when referred to a diagonal basis whose first $m$ members are invariant.

On the other hand, this projection takes $f_w(x)$ to $\frac 1{|G|}\sum e^{2\pi i\p w{Bx+b}}$ over
representatives $x\mapsto Bx+b$ of the point group $G\cong\G/\m{T}$. Since $\p w{Bx}=\p {B^{-1}w}x$ this is equally
$\frac 1{|G|}\sum e^{2\pi i \p wb} f_{B^{-1}w}(x)$, and so the $w,w$ entry in the matrix of $p$ will
be $\frac 1{|G|}\sum e^{2\pi i\p wb}$, summed over those $(B,b)$ for which $B^{-1}w=w$, or equivalently
$w=Bw$.

What this shows is that the multiplicity of $\mu$ in the Laplace spectrum of $M$ is
$\frac 1{|G|}\sum e^{2\pi i \p wb}$, summed over those pairs $w\in\m{T}^*$, $x\mapsto Bx+b$ for
which  $Bw=w$ and $4\pi^2N(w)=\mu$.

Using this formula, and collecting the terms with a given $\mu$, the trace of the heat kernel takes the form
\begin{eqnarray*}
\sum e^{-\mu t}&=&
\frac 1{|G|} \sum_{(B,b)} \,\,\, \sum_{w\in\m{T}^*:Bw=w} e^{2\pi i\p wb} e^{-4\pi^2 N(w) t} \\
&=&\frac 1{|G|} \sum_{(B,b)} \,\,\, \sum_{w\in{(p_B(\m{T}))}^*} g_t(w) \qquad\text{say,}
\end{eqnarray*}
since it is easily seen that the dual\ft{here $p_B$ denotes the orthogonal projection on the space fixed by $B$, and
the dual is taken in this space.} of $p_B(\m{T})$ is precisely the set of $w\in\m{T}^*$ such that $Bw=w$.
Now we use the Jacobi-Poisson formula that relates the values a Schwartz function $g$ takes in a lattice
to the values its Fourier transform $\widehat g$ takes in the dual lattice:
$$\sum_{x\in T} g(x) = \frac 1{\mathrm{vol}(\m{T})} \sum_{y\in \m{T}^*} \widehat{g}(y)$$
by applying it to each of the $|G|$ lattices ${(p_B(\m{T}))}^*$, getting
$$\frac 1{|G|} \sum_{(B,b)} \mathrm{vol}(p_B(\m{T})) \sum_{v\in p_B(\m{T})} \widehat{g}_t(v).$$
Then using the following rules for Fourier transforms of functions on $\R^m$:
\begin{itemize}
\item[] if $f(x)=e^{-\pi\p xx}$, then $\widehat{f}(y)=f(y)$,
\item[] if $g_\ld(x)=g(\ld x)$, $\ld>0$, then $\widehat{g_\ld}(y)=\frac 1{\ld^m}\widehat{g}(\frac y{\ld})$,
\item[] if $h_b(x)=e^{-2\pi i\p xb} h(x)$, then $\widehat{h_b}(y)=\widehat{h}(y-b)$,

\end{itemize}
we see first that the function $e^{-4\pi^2 tN(w)}$ on $\R^{n_\g}$ has Fourier transform
\newline
$(4\pi t)^{-\frac{n_\g}2} e^{-\frac 1{4t}\p vv}$ and second that its product with
$e^{2\pi i \p wb}=e^{2\pi i \p w{b_+}}$ has Fourier transform
$(4\pi t)^{-\frac{n_\g}2} e^{-\frac 1{4t}N(v+b_+)}$, where $b_+=p_B(b)$.\ft{Definitions of $\g$, $n_\g$ and $\m{T}_\g$ are
in Section \ref{hear}.}

The sum of this over all $v\in p_B(\m{T})$ is
$(4\pi t)^{-\frac{n_\g}2} \theta_{\m{T}_\g}(e^{-\frac 1{4t}})$,
giving the desired formula (\ref{trace}) for the trace of the heat kernel
$$\sum e^{-\mu t} = \frac 1{|G|} \sum_{\g\in G} \frac {\mathrm{vol}(\m{T}_\g)} {(4 \pi t)^{n_{\g}/2}}
\, \theta_{{\m T}_{\g}}(e^{\frac{-1}{4t}}).$$
This treatment is adapted from \cite{MRl}. Similar formulae appear in \cite{Su} and \cite{Gu}.

\

\noindent {\bf Note:} Isangulov \cite{Is} put an independent solution of the platycosm
isospectrality problem into the arXiv shortly after ours.
Although his statement is correct, the proof is incomplete,
the important error being the omission of the `balanced' vs `unbalanced'
cases for amphicosms that occupy roughly two-thirds of our argument.

\end{document}

%% file: fig1.pstex_t
\begin{picture}(0,0)%
\includegraphics{fig1.pstex}%
\end{picture}%
\setlength{\unitlength}{3947sp}%
\begingroup\makeatletter\ifx\SetFigFont\undefined%
\gdef\SetFigFont#1#2#3#4#5{%
  \reset@font\fontsize{#1}{#2pt}%
  \fontfamily{#3}\fontseries{#4}\fontshape{#5}%
  \selectfont}%
\fi\endgroup%
\begin{picture}(1797,1759)(1876,-1868)
\put(1876,-1814){\makebox(0,0)[lb]{\smash{\SetFigFont{12}{14.4}{\rmdefault}{\mddefault}{\updefault}$p_{01}$}}}
\put(3673,-1814){\makebox(0,0)[lb]{\smash{\SetFigFont{12}{14.4}{\rmdefault}{\mddefault}{\updefault}$p_{03}$}}}
\put(2774,-1814){\makebox(0,0)[lb]{\smash{\SetFigFont{12}{14.4}{\rmdefault}{\mddefault}{\updefault}$p_{13}$}}}
\put(2774,-253){\makebox(0,0)[lb]{\smash{\SetFigFont{12}{14.4}{\rmdefault}{\mddefault}{\updefault}$p_{02}$}}}
\put(2161,-963){\makebox(0,0)[lb]{\smash{\SetFigFont{12}{14.4}{\rmdefault}{\mddefault}{\updefault}$p_{12}$}}}
\put(3341,-963){\makebox(0,0)[lb]{\smash{\SetFigFont{12}{14.4}{\rmdefault}{\mddefault}{\updefault}$p_{23}$}}}
\put(2754,-1154){\makebox(0,0)[lb]{\smash{\SetFigFont{12}{14.4}{\rmdefault}{\mddefault}{\updefault}$0$}}}
\end{picture}

%% file: fig0.pstex_t
\begin{picture}(0,0)%
\includegraphics{fig0.pstex}%
\end{picture}%
\setlength{\unitlength}{3710sp}%
\begingroup\makeatletter\ifx\SetFigFont\undefined%
\gdef\SetFigFont#1#2#3#4#5{%
  \reset@font\fontsize{#1}{#2pt}%
  \fontfamily{#3}\fontseries{#4}\fontshape{#5}%
  \selectfont}%
\fi\endgroup%
\begin{picture}(6186,3427)(2195,-2868)
\put(2195,-1048){\makebox(0,0)[lb]{\smash{{\SetFigFont{9}{10.8}{\rmdefault}{\mddefault}{\updefault}$4$}}}}
\put(3861,-1048){\makebox(0,0)[lb]{\smash{{\SetFigFont{9}{10.8}{\rmdefault}{\mddefault}{\updefault}$6$}}}}
\put(2985,403){\makebox(0,0)[lb]{\smash{{\SetFigFont{9}{10.8}{\rmdefault}{\mddefault}{\updefault}$5$}}}}
\put(2951,-395){\makebox(0,0)[lb]{\smash{{\SetFigFont{9}{10.8}{\rmdefault}{\mddefault}{\updefault}$0$}}}}
\put(2966,-1048){\makebox(0,0)[lb]{\smash{{\SetFigFont{9}{10.8}{\rmdefault}{\mddefault}{\updefault}$-1$}}}}
\put(2463,-251){\makebox(0,0)[lb]{\smash{{\SetFigFont{9}{10.8}{\rmdefault}{\mddefault}{\updefault}$-1$}}}}
\put(3468,-251){\makebox(0,0)[lb]{\smash{{\SetFigFont{9}{10.8}{\rmdefault}{\mddefault}{\updefault}$-1$}}}}
\put(2391,324){\makebox(0,0)[lb]{\smash{{\SetFigFont{10}{12.0}{\rmdefault}{\mddefault}{\updefault}(a)}}}}
\put(6126,-1365){\makebox(0,0)[lb]{\smash{{\SetFigFont{9}{10.8}{\rmdefault}{\mddefault}{\updefault}$3$}}}}
\put(6937,-2817){\makebox(0,0)[lb]{\smash{{\SetFigFont{9}{10.8}{\rmdefault}{\mddefault}{\updefault}$2$}}}}
\put(5259,-2817){\makebox(0,0)[lb]{\smash{{\SetFigFont{9}{10.8}{\rmdefault}{\mddefault}{\updefault}$2$}}}}
\put(6054,-2817){\makebox(0,0)[lb]{\smash{{\SetFigFont{9}{10.8}{\rmdefault}{\mddefault}{\updefault}$1$}}}}
\put(6049,-2118){\makebox(0,0)[lb]{\smash{{\SetFigFont{9}{10.8}{\rmdefault}{\mddefault}{\updefault}$0$}}}}
\put(5546,-2046){\makebox(0,0)[lb]{\smash{{\SetFigFont{9}{10.8}{\rmdefault}{\mddefault}{\updefault}$-1$}}}}
\put(6571,-2046){\makebox(0,0)[lb]{\smash{{\SetFigFont{9}{10.8}{\rmdefault}{\mddefault}{\updefault}$1$}}}}
\put(5944,-1048){\makebox(0,0)[lb]{\smash{{\SetFigFont{9}{10.8}{\rmdefault}{\mddefault}{\updefault}$3$}}}}
\put(5575,-251){\makebox(0,0)[lb]{\smash{{\SetFigFont{9}{10.8}{\rmdefault}{\mddefault}{\updefault}$0$}}}}
\put(4209,-1048){\makebox(0,0)[lb]{\smash{{\SetFigFont{9}{10.8}{\rmdefault}{\mddefault}{\updefault}$1$}}}}
\put(5031,-1048){\makebox(0,0)[lb]{\smash{{\SetFigFont{9}{10.8}{\rmdefault}{\mddefault}{\updefault}$2$}}}}
\put(5109,395){\makebox(0,0)[lb]{\smash{{\SetFigFont{9}{10.8}{\rmdefault}{\mddefault}{\updefault}$4$}}}}
\put(4956,-371){\makebox(0,0)[lb]{\smash{{\SetFigFont{9}{10.8}{\rmdefault}{\mddefault}{\updefault}$-1$}}}}
\put(4616,-251){\makebox(0,0)[lb]{\smash{{\SetFigFont{9}{10.8}{\rmdefault}{\mddefault}{\updefault}$0$}}}}
\put(4381,324){\makebox(0,0)[lb]{\smash{{\SetFigFont{10}{12.0}{\rmdefault}{\mddefault}{\updefault}(c)}}}}
\put(5176,-2386){\makebox(0,0)[lb]{\smash{{\SetFigFont{10}{12.0}{\rmdefault}{\mddefault}{\updefault}(d)}}}}
\put(4823,-2828){\makebox(0,0)[lb]{\smash{{\SetFigFont{9}{10.8}{\rmdefault}{\mddefault}{\updefault}$5$}}}}
\put(3133,-2828){\makebox(0,0)[lb]{\smash{{\SetFigFont{9}{10.8}{\rmdefault}{\mddefault}{\updefault}$3$}}}}
\put(3995,-2828){\makebox(0,0)[lb]{\smash{{\SetFigFont{9}{10.8}{\rmdefault}{\mddefault}{\updefault}$0$}}}}
\put(4449,-2031){\makebox(0,0)[lb]{\smash{{\SetFigFont{9}{10.8}{\rmdefault}{\mddefault}{\updefault}$-2$}}}}
\put(3451,-2031){\makebox(0,0)[lb]{\smash{{\SetFigFont{9}{10.8}{\rmdefault}{\mddefault}{\updefault}$-2$}}}}
\put(3923,-2175){\makebox(0,0)[lb]{\smash{{\SetFigFont{9}{10.8}{\rmdefault}{\mddefault}{\updefault}$1$}}}}
\put(3995,-1385){\makebox(0,0)[lb]{\smash{{\SetFigFont{9}{10.8}{\rmdefault}{\mddefault}{\updefault}$6$}}}}
\put(2926,-2441){\makebox(0,0)[lb]{\smash{{\SetFigFont{10}{12.0}{\rmdefault}{\mddefault}{\updefault}(b)}}}}
\put(7323,-1057){\makebox(0,0)[lb]{\smash{{\SetFigFont{9}{10.8}{\rmdefault}{\mddefault}{\updefault}$0$}}}}
\put(6488,-1057){\makebox(0,0)[lb]{\smash{{\SetFigFont{9}{10.8}{\rmdefault}{\mddefault}{\updefault}$1$}}}}
\put(8151,-1057){\makebox(0,0)[lb]{\smash{{\SetFigFont{9}{10.8}{\rmdefault}{\mddefault}{\updefault}$3$}}}}
\put(7309,451){\makebox(0,0)[lb]{\smash{{\SetFigFont{9}{10.8}{\rmdefault}{\mddefault}{\updefault}$2$}}}}
\put(7228,-349){\makebox(0,0)[lb]{\smash{{\SetFigFont{9}{10.8}{\rmdefault}{\mddefault}{\updefault}$1$}}}}
\put(7788,-276){\makebox(0,0)[lb]{\smash{{\SetFigFont{9}{10.8}{\rmdefault}{\mddefault}{\updefault}$0$}}}}
\put(6807,-276){\makebox(0,0)[lb]{\smash{{\SetFigFont{9}{10.8}{\rmdefault}{\mddefault}{\updefault}$0$}}}}
\put(6653,305){\makebox(0,0)[lb]{\smash{{\SetFigFont{10}{12.0}{\rmdefault}{\mddefault}{\updefault}(e)}}}}
\end{picture}%

%% file: fig6.pstex_t
\begin{picture}(0,0)%
\includegraphics{fig6.pstex}%
\end{picture}%
\setlength{\unitlength}{3947sp}%
\begingroup\makeatletter\ifx\SetFigFont\undefined%
\gdef\SetFigFont#1#2#3#4#5{%
  \reset@font\fontsize{#1}{#2pt}%
  \fontfamily{#3}\fontseries{#4}\fontshape{#5}%
  \selectfont}%
\fi\endgroup%
\begin{picture}(1797,1696)(1876,-1874)
\put(1876,-1814){\makebox(0,0)[lb]{\smash{\SetFigFont{11}{13.2}{\rmdefault}{\mddefault}{\updefault}$A$}}}
\put(2691,-1101){\makebox(0,0)[lb]{\smash{\SetFigFont{11}{13.2}{\rmdefault}{\mddefault}{\updefault}$E$}}}
\put(2774,-1834){\makebox(0,0)[lb]{\smash{\SetFigFont{11}{13.2}{\rmdefault}{\mddefault}{\updefault}$B$}}}
\put(3673,-1814){\makebox(0,0)[lb]{\smash{\SetFigFont{11}{13.2}{\rmdefault}{\mddefault}{\updefault}$C$}}}
\put(3281,-981){\makebox(0,0)[lb]{\smash{\SetFigFont{11}{13.2}{\rmdefault}{\mddefault}{\updefault}$F$}}}
\put(2796,-286){\makebox(0,0)[lb]{\smash{\SetFigFont{11}{13.2}{\rmdefault}{\mddefault}{\updefault}$0$}}}
\put(2281,-983){\makebox(0,0)[lb]{\smash{\SetFigFont{11}{13.2}{\rmdefault}{\mddefault}{\updefault}$D$}}}
\end{picture}

%% file: fig11.pstex_t
\begin{picture}(0,0)%
\includegraphics{fig11.pstex}%
\end{picture}%
\setlength{\unitlength}{3947sp}%
\begingroup\makeatletter\ifx\SetFigFont\undefined%
\gdef\SetFigFont#1#2#3#4#5{%
  \reset@font\fontsize{#1}{#2pt}%
  \fontfamily{#3}\fontseries{#4}\fontshape{#5}%
  \selectfont}%
\fi\endgroup%
\begin{picture}(5768,2011)(751,-1986)
\put(751,-1061){\makebox(0,0)[lb]{\smash{{\SetFigFont{14}{16.8}{\rmdefault}{\mddefault}{\updefault}$\x$}}}}
\put(1573,-130){\makebox(0,0)[lb]{\smash{{\SetFigFont{14}{16.8}{\rmdefault}{\mddefault}{\updefault}$\w$}}}}
\put(1403,-1939){\makebox(0,0)[lb]{\smash{{\SetFigFont{14}{16.8}{\rmdefault}{\mddefault}{\updefault}$\z$}}}}
\put(1548,-1939){\makebox(0,0)[lb]{\smash{{\SetFigFont{12}{14.4}{\rmdefault}{\mddefault}{\updefault}(norm $D$)}}}}
\put(4782,-103){\makebox(0,0)[lb]{\smash{{\SetFigFont{14}{16.8}{\rmdefault}{\mddefault}{\updefault}$\w$}}}}
\put(3905,-1061){\makebox(0,0)[lb]{\smash{{\SetFigFont{14}{16.8}{\rmdefault}{\mddefault}{\updefault}$\x$}}}}
\put(4589,-1939){\makebox(0,0)[lb]{\smash{{\SetFigFont{14}{16.8}{\rmdefault}{\mddefault}{\updefault}$\z$}}}}
\put(4734,-1939){\makebox(0,0)[lb]{\smash{{\SetFigFont{12}{14.4}{\rmdefault}{\mddefault}{\updefault}(norm $D$)}}}}
\end{picture}%

%% file: fig12.pstex_t
\begin{picture}(0,0)%
\includegraphics{fig12.pstex}%
\end{picture}%
\setlength{\unitlength}{3947sp}%
\begingroup\makeatletter\ifx\SetFigFont\undefined%
\gdef\SetFigFont#1#2#3#4#5{%
  \reset@font\fontsize{#1}{#2pt}%
  \fontfamily{#3}\fontseries{#4}\fontshape{#5}%
  \selectfont}%
\fi\endgroup%
\begin{picture}(2112,1789)(1126,-1844)
\put(2176,-211){\makebox(0,0)[lb]{\smash{\SetFigFont{12}{14.4}{\rmdefault}{\mddefault}{\updefault}${\m D}=\Ld_{A\,B\,C}$}}}
\put(3151,-961){\makebox(0,0)[lb]{\smash{\SetFigFont{12}{14.4}{\familydefault}{\mddefault}{\updefault}$\Ld_{A\,B\,[C]}$}}}
\put(2326,-961){\makebox(0,0)[lb]{\smash{\SetFigFont{12}{14.4}{\familydefault}{\mddefault}{\updefault}$\Ld_{A\,[B]\,C}$}}}
\put(1126,-961){\makebox(0,0)[lb]{\smash{\SetFigFont{12}{14.4}{\familydefault}{\mddefault}{\updefault}${\m H}=\Ld_{[A]\,B\,C}$}}}
\put(2326,-1786){\makebox(0,0)[lb]{\smash{\SetFigFont{12}{14.4}{\familydefault}{\mddefault}{\updefault}$\Ld_{4A\,4B\,4C}$}}}
\end{picture}

%% file: fig13.pstex_t
\begin{picture}(0,0)%
\includegraphics{fig13.pstex}%
\end{picture}%
\setlength{\unitlength}{3947sp}%
\begingroup\makeatletter\ifx\SetFigFont\undefined%
\gdef\SetFigFont#1#2#3#4#5{%
  \reset@font\fontsize{#1}{#2pt}%
  \fontfamily{#3}\fontseries{#4}\fontshape{#5}%
  \selectfont}%
\fi\endgroup%
\begin{picture}(5954,3334)(1233,-2672)
\put(4102,-228){\makebox(0,0)[lb]{\smash{{\SetFigFont{9}{10.8}{\rmdefault}{\mddefault}{\updefault}$0$}}}}
\put(3177,-785){\makebox(0,0)[lb]{\smash{{\SetFigFont{9}{10.8}{\rmdefault}{\mddefault}{\updefault}$D\!+\!B\!-\!C$}}}}
\put(3866,-855){\makebox(0,0)[lb]{\smash{{\SetFigFont{9}{10.8}{\rmdefault}{\mddefault}{\updefault}$B\!+\!C\!\!-\!\!D$}}}}
\put(3036,-1535){\makebox(0,0)[lb]{\smash{{\SetFigFont{9}{10.8}{\rmdefault}{\mddefault}{\updefault}$D\!+\!B\!-\!C$}}}}
\put(4730,-1535){\makebox(0,0)[lb]{\smash{{\SetFigFont{9}{10.8}{\rmdefault}{\mddefault}{\updefault}$D\!+\!C\!-\!B$}}}}
\put(3742,-1579){\makebox(0,0)[lb]{\smash{{\SetFigFont{9}{10.8}{\rmdefault}{\mddefault}{\updefault}$4A\!+\!B\!+\!C\!\!-\!D$}}}}
\put(2727,-1535){\makebox(0,0)[lb]{\smash{{\SetFigFont{9}{10.8}{\rmdefault}{\mddefault}{\updefault}$2C$}}}}
\put(1945,-963){\makebox(0,0)[lb]{\smash{{\SetFigFont{9}{10.8}{\rmdefault}{\mddefault}{\updefault}$0$}}}}
\put(1713,-228){\makebox(0,0)[lb]{\smash{{\SetFigFont{9}{10.8}{\rmdefault}{\mddefault}{\updefault}$D\!-\!B\!-\!C$}}}}
\put(1233,-1535){\makebox(0,0)[lb]{\smash{{\SetFigFont{9}{10.8}{\rmdefault}{\mddefault}{\updefault}$2B$}}}}
\put(1492,-823){\makebox(0,0)[lb]{\smash{{\SetFigFont{9}{10.8}{\rmdefault}{\mddefault}{\updefault}$2B$}}}}
\put(2410,-829){\makebox(0,0)[lb]{\smash{{\SetFigFont{9}{10.8}{\rmdefault}{\mddefault}{\updefault}$2C$}}}}
\put(1951,-1579){\makebox(0,0)[lb]{\smash{{\SetFigFont{9}{10.8}{\rmdefault}{\mddefault}{\updefault}$4A$}}}}
\put(5374,-740){\makebox(0,0)[lb]{\smash{{\SetFigFont{9}{10.8}{\rmdefault}{\mddefault}{\updefault}$2D$}}}}
\put(5930,566){\makebox(0,0)[lb]{\smash{{\SetFigFont{9}{10.8}{\rmdefault}{\mddefault}{\updefault}$B\!-\!C\!-\!D$}}}}
\put(6023,-105){\makebox(0,0)[lb]{\smash{{\SetFigFont{9}{10.8}{\rmdefault}{\mddefault}{\updefault}$2C$}}}}
\put(6889,-740){\makebox(0,0)[lb]{\smash{{\SetFigFont{9}{10.8}{\rmdefault}{\mddefault}{\updefault}$0$}}}}
\put(5966,-785){\makebox(0,0)[lb]{\smash{{\SetFigFont{9}{10.8}{\rmdefault}{\mddefault}{\updefault}$4A\!+\!2C$}}}}
\put(5436,-2593){\makebox(0,0)[lb]{\smash{{\SetFigFont{9}{10.8}{\rmdefault}{\mddefault}{\updefault}$0$}}}}
\put(5763,-1900){\makebox(0,0)[lb]{\smash{{\SetFigFont{9}{10.8}{\rmdefault}{\mddefault}{\updefault}$0$}}}}
\put(6581,-1900){\makebox(0,0)[lb]{\smash{{\SetFigFont{9}{10.8}{\rmdefault}{\mddefault}{\updefault}$2D$}}}}
\put(6023,-1958){\makebox(0,0)[lb]{\smash{{\SetFigFont{9}{10.8}{\rmdefault}{\mddefault}{\updefault}$2B$}}}}
\put(6889,-2593){\makebox(0,0)[lb]{\smash{{\SetFigFont{9}{10.8}{\rmdefault}{\mddefault}{\updefault}$2D$}}}}
\put(5966,-2637){\makebox(0,0)[lb]{\smash{{\SetFigFont{9}{10.8}{\rmdefault}{\mddefault}{\updefault}$4A\!+\!2B$}}}}
\put(5930,-1287){\makebox(0,0)[lb]{\smash{{\SetFigFont{9}{10.8}{\rmdefault}{\mddefault}{\updefault}$C\!-\!B\!-\!D$}}}}
\put(6601,-47){\makebox(0,0)[lb]{\smash{{\SetFigFont{9}{10.8}{\rmdefault}{\mddefault}{\updefault}$0$}}}}
\put(4527,-785){\makebox(0,0)[lb]{\smash{{\SetFigFont{9}{10.8}{\rmdefault}{\mddefault}{\updefault}$D\!+\!C\!-\!B$}}}}
\put(5700,-47){\makebox(0,0)[lb]{\smash{{\SetFigFont{9}{10.8}{\rmdefault}{\mddefault}{\updefault}$2D$}}}}
\end{picture}%

%% file: fig14.pstex_t
\begin{picture}(0,0)%
\includegraphics{fig14.pstex}%
\end{picture}%
\setlength{\unitlength}{3947sp}%
\begingroup\makeatletter\ifx\SetFigFont\undefined%
\gdef\SetFigFont#1#2#3#4#5{%
  \reset@font\fontsize{#1}{#2pt}%
  \fontfamily{#3}\fontseries{#4}\fontshape{#5}%
  \selectfont}%
\fi\endgroup%
\begin{picture}(1947,1588)(1351,-1076)
\put(2476,-1036){\makebox(0,0)[lb]{\smash{\SetFigFont{10}{12.0}{\rmdefault}{\mddefault}{\updefault}$2B$}}}
\put(3286,-1036){\makebox(0,0)[lb]{\smash{\SetFigFont{10}{12.0}{\rmdefault}{\mddefault}{\updefault}$C-B$}}}
\put(3016,-271){\makebox(0,0)[lb]{\smash{\SetFigFont{10}{12.0}{\rmdefault}{\mddefault}{\updefault}$0$}}}
\put(2476,404){\makebox(0,0)[lb]{\smash{\SetFigFont{10}{12.0}{\rmdefault}{\mddefault}{\updefault}$4D$}}}
\put(2071,-271){\makebox(0,0)[lb]{\smash{\SetFigFont{10}{12.0}{\rmdefault}{\mddefault}{\updefault}$0$}}}
\put(2470,-406){\makebox(0,0)[lb]{\smash{\SetFigFont{10}{12.0}{\rmdefault}{\mddefault}{\updefault}$0$}}}
\put(1351,-1036){\makebox(0,0)[lb]{\smash{\SetFigFont{10}{12.0}{\rmdefault}{\mddefault}{\updefault}$4A+2B$}}}
\end{picture}

%% file: fig15.pstex_t
\begin{picture}(0,0)%
\includegraphics{fig15.pstex}%
\end{picture}%
\setlength{\unitlength}{3947sp}%
\begingroup\makeatletter\ifx\SetFigFont\undefined%
\gdef\SetFigFont#1#2#3#4#5{%
  \reset@font\fontsize{#1}{#2pt}%
  \fontfamily{#3}\fontseries{#4}\fontshape{#5}%
  \selectfont}%
\fi\endgroup%
\begin{picture}(5938,1735)(1576,-1247)
\put(4407,242){\makebox(0,0)[lb]{\smash{{\SetFigFont{11}{13.2}{\rmdefault}{\mddefault}{\updefault}$0$}}}}
\put(3415,-356){\makebox(0,0)[lb]{\smash{{\SetFigFont{10}{12.0}{\rmdefault}{\mddefault}{\updefault}$D\!+\!B\!-\!C$}}}}
\put(4155,-432){\makebox(0,0)[lb]{\smash{{\SetFigFont{10}{12.0}{\rmdefault}{\mddefault}{\updefault}$B\!+\!C\!-\!\!D$}}}}
\put(3264,-1161){\makebox(0,0)[lb]{\smash{{\SetFigFont{10}{12.0}{\rmdefault}{\mddefault}{\updefault}$D\!+\!B\!-\!C$}}}}
\put(5084,-1161){\makebox(0,0)[lb]{\smash{{\SetFigFont{10}{12.0}{\rmdefault}{\mddefault}{\updefault}$D\!+\!C\!-\!B$}}}}
\put(4848,-356){\makebox(0,0)[lb]{\smash{{\SetFigFont{10}{12.0}{\rmdefault}{\mddefault}{\updefault}$D\!+\!C\!-\!B$}}}}
\put(4022,-1209){\makebox(0,0)[lb]{\smash{{\SetFigFont{10}{12.0}{\rmdefault}{\mddefault}{\updefault}$4A\!+\!B\!+\!C\!\!-\!D$}}}}
\put(2304,-405){\makebox(0,0)[lb]{\smash{{\SetFigFont{11}{13.2}{\rmdefault}{\mddefault}{\updefault}$0$}}}}
\put(2055,384){\makebox(0,0)[lb]{\smash{{\SetFigFont{10}{12.0}{\rmdefault}{\mddefault}{\updefault}$D\!-\!B\!-\!C$}}}}
\put(1818,-255){\makebox(0,0)[lb]{\smash{{\SetFigFont{10}{12.0}{\rmdefault}{\mddefault}{\updefault}$2B$}}}}
\put(2804,-261){\makebox(0,0)[lb]{\smash{{\SetFigFont{10}{12.0}{\rmdefault}{\mddefault}{\updefault}$2C$}}}}
\put(2311,-1067){\makebox(0,0)[lb]{\smash{{\SetFigFont{10}{12.0}{\rmdefault}{\mddefault}{\updefault}$4A$}}}}
\put(2998,-1048){\makebox(0,0)[lb]{\smash{{\SetFigFont{10}{12.0}{\rmdefault}{\mddefault}{\updefault}$2C$}}}}
\put(6052,-275){\makebox(0,0)[lb]{\smash{{\SetFigFont{11}{13.2}{\rmdefault}{\mddefault}{\updefault}$0$}}}}
\put(6931,-275){\makebox(0,0)[lb]{\smash{{\SetFigFont{10}{12.0}{\rmdefault}{\mddefault}{\updefault}$2D$}}}}
\put(6231,384){\makebox(0,0)[lb]{\smash{{\SetFigFont{10}{12.0}{\rmdefault}{\mddefault}{\updefault}$C\!-\!B\!-\!D$}}}}
\put(6329,-337){\makebox(0,0)[lb]{\smash{{\SetFigFont{10}{12.0}{\rmdefault}{\mddefault}{\updefault}$2B$}}}}
\put(6269,-1067){\makebox(0,0)[lb]{\smash{{\SetFigFont{10}{12.0}{\rmdefault}{\mddefault}{\updefault}$4A\!+\!2B$}}}}
\put(1576,-1048){\makebox(0,0)[lb]{\smash{{\SetFigFont{10}{12.0}{\rmdefault}{\mddefault}{\updefault}$2B$}}}}
\put(7193,-1048){\makebox(0,0)[lb]{\smash{{\SetFigFont{10}{12.0}{\rmdefault}{\mddefault}{\updefault}$2D$}}}}
\put(5771,-1029){\makebox(0,0)[lb]{\smash{{\SetFigFont{11}{13.2}{\rmdefault}{\mddefault}{\updefault}$0$}}}}
\end{picture}%

%% file: fig16.pstex_t
\begin{picture}(0,0)%
\includegraphics{fig16.pstex}%
\end{picture}%
\setlength{\unitlength}{3947sp}%
\begingroup\makeatletter\ifx\SetFigFont\undefined%
\gdef\SetFigFont#1#2#3#4#5{%
  \reset@font\fontsize{#1}{#2pt}%
  \fontfamily{#3}\fontseries{#4}\fontshape{#5}%
  \selectfont}%
\fi\endgroup%
\begin{picture}(4258,1603)(1726,-1993)
\put(3597,-569){\makebox(0,0)[lb]{\smash{\SetFigFont{14}{16.8}{\rmdefault}{\mddefault}{\updefault}$\m{T}\ne\m{T}'$}}}
\put(1726,-1953){\makebox(0,0)[lb]{\smash{\SetFigFont{10}{12.0}{\rmdefault}{\mddefault}{\updefault}$4A$}}}
\put(2614,-1953){\makebox(0,0)[lb]{\smash{\SetFigFont{10}{12.0}{\rmdefault}{\mddefault}{\updefault}$0$}}}
\put(3390,-1953){\makebox(0,0)[lb]{\smash{\SetFigFont{10}{12.0}{\rmdefault}{\mddefault}{\updefault}$C$}}}
\put(3115,-1179){\makebox(0,0)[lb]{\smash{\SetFigFont{10}{12.0}{\rmdefault}{\mddefault}{\updefault}$0$}}}
\put(2569,-498){\makebox(0,0)[lb]{\smash{\SetFigFont{10}{12.0}{\rmdefault}{\mddefault}{\updefault}$4C$}}}
\put(2160,-1179){\makebox(0,0)[lb]{\smash{\SetFigFont{10}{12.0}{\rmdefault}{\mddefault}{\updefault}$0$}}}
\put(4336,-1953){\makebox(0,0)[lb]{\smash{\SetFigFont{10}{12.0}{\rmdefault}{\mddefault}{\updefault}$0$}}}
\put(5119,-1953){\makebox(0,0)[lb]{\smash{\SetFigFont{10}{12.0}{\rmdefault}{\mddefault}{\updefault}$4A$}}}
\put(5984,-1953){\makebox(0,0)[lb]{\smash{\SetFigFont{10}{12.0}{\rmdefault}{\mddefault}{\updefault}$2C$}}}
\put(5627,-1192){\makebox(0,0)[lb]{\smash{\SetFigFont{10}{12.0}{\rmdefault}{\mddefault}{\updefault}$2C$}}}
\put(4656,-1179){\makebox(0,0)[lb]{\smash{\SetFigFont{10}{12.0}{\rmdefault}{\mddefault}{\updefault}$0$}}}
\put(5166,-498){\makebox(0,0)[lb]{\smash{\SetFigFont{10}{12.0}{\rmdefault}{\mddefault}{\updefault}$0$}}}
\put(2563,-1317){\makebox(0,0)[lb]{\smash{\SetFigFont{10}{12.0}{\rmdefault}{\mddefault}{\updefault}$0$}}}
\put(5113,-1317){\makebox(0,0)[lb]{\smash{\SetFigFont{10}{12.0}{\rmdefault}{\mddefault}{\updefault}$0$}}}
\end{picture}

%% file: fig17.pstex_t
\begin{picture}(0,0)%
\includegraphics{fig17.pstex}%
\end{picture}%
\setlength{\unitlength}{3947sp}%
\begingroup\makeatletter\ifx\SetFigFont\undefined%
\gdef\SetFigFont#1#2#3#4#5{%
  \reset@font\fontsize{#1}{#2pt}%
  \fontfamily{#3}\fontseries{#4}\fontshape{#5}%
  \selectfont}%
\fi\endgroup%
\begin{picture}(3630,1834)(733,-1873)
\put(2926,-886){\makebox(0,0)[lb]{\smash{\SetFigFont{14}{16.8}{\rmdefault}{\mddefault}{\updefault}$\vv$}}}
\put(2926,-1786){\makebox(0,0)[lb]{\smash{\SetFigFont{14}{16.8}{\rmdefault}{\mddefault}{\updefault}$\w$}}}
\put(1876,-886){\makebox(0,0)[lb]{\smash{\SetFigFont{14}{16.8}{\rmdefault}{\mddefault}{\updefault}$-\w$}}}
\put(1876,-1786){\makebox(0,0)[lb]{\smash{\SetFigFont{14}{16.8}{\rmdefault}{\mddefault}{\updefault}$-\vv$}}}
\put(2531,-1336){\makebox(0,0)[lb]{\smash{\SetFigFont{14}{16.8}{\rmdefault}{\mddefault}{\updefault}$0$}}}
\end{picture}

%% file: fig18.pstex_t
\begin{picture}(0,0)%
\includegraphics{fig18.pstex}%
\end{picture}%
\setlength{\unitlength}{3947sp}%
\begingroup\makeatletter\ifx\SetFigFont\undefined%
\gdef\SetFigFont#1#2#3#4#5{%
  \reset@font\fontsize{#1}{#2pt}%
  \fontfamily{#3}\fontseries{#4}\fontshape{#5}%
  \selectfont}%
\fi\endgroup%
\begin{picture}(5187,1665)(1201,-1756)
\put(3526,-1711){\makebox(0,0)[lb]{\smash{{\SetFigFont{12}{14.4}{\familydefault}{\mddefault}{\updefault}$(a)$}}}}
\put(3151,-961){\makebox(0,0)[lb]{\smash{{\SetFigFont{12}{14.4}{\familydefault}{\mddefault}{\updefault}$\Ld_{\al\,\be\,[\be]}$}}}}
\put(2326,-961){\makebox(0,0)[lb]{\smash{{\SetFigFont{12}{14.4}{\familydefault}{\mddefault}{\updefault}$\Ld_{\al\,[\be]\,\be}$}}}}
\put(1201,-961){\makebox(0,0)[lb]{\smash{{\SetFigFont{12}{14.4}{\familydefault}{\mddefault}{\updefault}${\m H}=\Ld_{[\al]\,\be\,\be}$}}}}
\put(2101,-211){\makebox(0,0)[lb]{\smash{{\SetFigFont{12}{14.4}{\rmdefault}{\mddefault}{\updefault}${\m D}=\Ld_{\al\,\be\,\be}$}}}}
\put(2326,-1711){\makebox(0,0)[lb]{\smash{{\SetFigFont{12}{14.4}{\familydefault}{\mddefault}{\updefault}$\Ld_{4\al\,4\be\,4\be}$}}}}
\put(4501,-1711){\makebox(0,0)[lb]{\smash{{\SetFigFont{12}{14.4}{\familydefault}{\mddefault}{\updefault}$(b)$}}}}
\put(4051,-511){\makebox(0,0)[lb]{\smash{{\SetFigFont{12}{14.4}{\familydefault}{\mddefault}{\updefault}$\te_2+\te_3=2\te_3$}}}}
\put(4501,-929){\makebox(0,0)[lb]{\smash{{\SetFigFont{12}{14.4}{\familydefault}{\mddefault}{\updefault}$\te_0+\te_1$}}}}
\put(5326,-929){\makebox(0,0)[lb]{\smash{{\SetFigFont{12}{14.4}{\familydefault}{\mddefault}{\updefault}$\te_0+\te_2$}}}}
\put(5551,-1691){\makebox(0,0)[lb]{\smash{{\SetFigFont{12}{14.4}{\familydefault}{\mddefault}{\updefault}$\te_0$}}}}
\put(5026,-211){\makebox(0,0)[lb]{\smash{{\SetFigFont{12}{14.4}{\familydefault}{\mddefault}{\updefault}$\te_0+\te_1+\te_2+\te_3$}}}}
\put(5661,-1290){\makebox(0,0)[lb]{\smash{{\SetFigFont{12}{14.4}{\rmdefault}{\mddefault}{\updefault}$=$}}}}
\put(5851,-1336){\makebox(0,0)[lb]{\smash{{\SetFigFont{12}{14.4}{\familydefault}{\mddefault}{\updefault}$\te_3$}}}}
\put(5401,-1261){\makebox(0,0)[lb]{\smash{{\SetFigFont{12}{14.4}{\familydefault}{\mddefault}{\updefault}$\te_2$}}}}
\put(6151,-961){\makebox(0,0)[lb]{\smash{{\SetFigFont{12}{14.4}{\familydefault}{\mddefault}{\updefault}$\te_0+\te_3$}}}}
\end{picture}%

%% file: fig19.pstex_t
\begin{picture}(0,0)%
\includegraphics{fig19.pstex}%
\end{picture}%
\setlength{\unitlength}{3947sp}%
\begingroup\makeatletter\ifx\SetFigFont\undefined%
\gdef\SetFigFont#1#2#3#4#5{%
  \reset@font\fontsize{#1}{#2pt}%
  \fontfamily{#3}\fontseries{#4}\fontshape{#5}%
  \selectfont}%
\fi\endgroup%
\begin{picture}(5698,1608)(1576,-1152)
\put(3603,-1060){\makebox(0,0)[lb]{\smash{{\SetFigFont{9}{10.8}{\rmdefault}{\mddefault}{\updefault}$2$}}}}
\put(5068,-1060){\makebox(0,0)[lb]{\smash{{\SetFigFont{9}{10.8}{\rmdefault}{\mddefault}{\updefault}$2$}}}}
\put(4273,279){\makebox(0,0)[lb]{\smash{{\SetFigFont{9}{10.8}{\rmdefault}{\mddefault}{\updefault}$4\delta-2$}}}}
\put(3938,-348){\makebox(0,0)[lb]{\smash{{\SetFigFont{9}{10.8}{\rmdefault}{\mddefault}{\updefault}$2$}}}}
\put(4776,-348){\makebox(0,0)[lb]{\smash{{\SetFigFont{9}{10.8}{\rmdefault}{\mddefault}{\updefault}$2$}}}}
\put(4309,-473){\makebox(0,0)[lb]{\smash{{\SetFigFont{9}{10.8}{\rmdefault}{\mddefault}{\updefault}$0$}}}}
\put(4315,-1060){\makebox(0,0)[lb]{\smash{{\SetFigFont{9}{10.8}{\rmdefault}{\mddefault}{\updefault}$4\al$}}}}
\put(3099,-37){\makebox(0,0)[lb]{\smash{{\SetFigFont{11}{13.2}{\familydefault}{\mddefault}{\updefault}{\bf or}}}}}
\put(5177,241){\makebox(0,0)[lb]{\smash{{\SetFigFont{11}{13.2}{\rmdefault}{\mddefault}{\updefault}${\m T}_-$}}}}
\put(2358,279){\makebox(0,0)[lb]{\smash{{\SetFigFont{9}{10.8}{\rmdefault}{\mddefault}{\updefault}$16\delta$}}}}
\put(2830,-360){\makebox(0,0)[lb]{\smash{{\SetFigFont{9}{10.8}{\rmdefault}{\mddefault}{\updefault}$0$}}}}
\put(2353,-473){\makebox(0,0)[lb]{\smash{{\SetFigFont{9}{10.8}{\rmdefault}{\mddefault}{\updefault}$0$}}}}
\put(1981,-348){\makebox(0,0)[lb]{\smash{{\SetFigFont{9}{10.8}{\rmdefault}{\mddefault}{\updefault}$0$}}}}
\put(2400,-1060){\makebox(0,0)[lb]{\smash{{\SetFigFont{9}{10.8}{\rmdefault}{\mddefault}{\updefault}$2$}}}}
\put(1645,241){\makebox(0,0)[lb]{\smash{{\SetFigFont{11}{13.2}{\rmdefault}{\mddefault}{\updefault}${\m T}_+$}}}}
\put(3165,-1060){\makebox(0,0)[lb]{\smash{{\SetFigFont{9}{10.8}{\rmdefault}{\mddefault}{\updefault}$0$}}}}
\put(1576,-1074){\makebox(0,0)[lb]{\smash{{\SetFigFont{9}{10.8}{\rmdefault}{\mddefault}{\updefault}$4\al+2$}}}}
\put(5923,-1060){\makebox(0,0)[lb]{\smash{{\SetFigFont{9}{10.8}{\rmdefault}{\mddefault}{\updefault}$2+4\al-4\delta$}}}}
\put(6139,279){\makebox(0,0)[lb]{\smash{{\SetFigFont{9}{10.8}{\rmdefault}{\mddefault}{\updefault}$0$}}}}
\put(5672,-348){\makebox(0,0)[lb]{\smash{{\SetFigFont{9}{10.8}{\rmdefault}{\mddefault}{\updefault}$4\delta$}}}}
\put(6008,-432){\makebox(0,0)[lb]{\smash{{\SetFigFont{9}{10.8}{\rmdefault}{\mddefault}{\updefault}$2-4\delta$}}}}
\put(6552,-348){\makebox(0,0)[lb]{\smash{{\SetFigFont{9}{10.8}{\rmdefault}{\mddefault}{\updefault}$4\delta$}}}}
\put(5374,-1060){\makebox(0,0)[lb]{\smash{{\SetFigFont{9}{10.8}{\rmdefault}{\mddefault}{\updefault}$4\delta$}}}}
\put(6893,-1060){\makebox(0,0)[lb]{\smash{{\SetFigFont{9}{10.8}{\rmdefault}{\mddefault}{\updefault}$4\delta$}}}}
\end{picture}%

%% file: fig20.pstex_t
\begin{picture}(0,0)%
\includegraphics{fig20.pstex}%
\end{picture}%
\setlength{\unitlength}{3947sp}%
\begingroup\makeatletter\ifx\SetFigFont\undefined%
\gdef\SetFigFont#1#2#3#4#5{%
  \reset@font\fontsize{#1}{#2pt}%
  \fontfamily{#3}\fontseries{#4}\fontshape{#5}%
  \selectfont}%
\fi\endgroup%
\begin{picture}(5718,2917)(1801,-2356)
\put(5315,-801){\makebox(0,0)[lb]{\smash{{\SetFigFont{9}{10.8}{\rmdefault}{\mddefault}{\updefault}$4$}}}}
\put(4996,-142){\makebox(0,0)[lb]{\smash{{\SetFigFont{9}{10.8}{\rmdefault}{\mddefault}{\updefault}$4$}}}}
\put(4146,-142){\makebox(0,0)[lb]{\smash{{\SetFigFont{9}{10.8}{\rmdefault}{\mddefault}{\updefault}$12$}}}}
\put(6851,-142){\makebox(0,0)[lb]{\smash{{\SetFigFont{9}{10.8}{\rmdefault}{\mddefault}{\updefault}$4\al$}}}}
\put(4614,-1071){\makebox(0,0)[lb]{\smash{{\SetFigFont{9}{10.8}{\rmdefault}{\mddefault}{\updefault}$0$}}}}
\put(6001,-1651){\makebox(0,0)[lb]{\smash{{\SetFigFont{9}{10.8}{\rmdefault}{\mddefault}{\updefault}$2\delta$}}}}
\put(6862,-1646){\makebox(0,0)[lb]{\smash{{\SetFigFont{9}{10.8}{\rmdefault}{\mddefault}{\updefault}$0$}}}}
\put(5754,-2309){\makebox(0,0)[lb]{\smash{{\SetFigFont{9}{10.8}{\rmdefault}{\mddefault}{\updefault}$2\delta$}}}}
\put(7171,-2309){\makebox(0,0)[lb]{\smash{{\SetFigFont{9}{10.8}{\rmdefault}{\mddefault}{\updefault}$0$}}}}
\put(4540,437){\makebox(0,0)[lb]{\smash{{\SetFigFont{9}{10.8}{\rmdefault}{\mddefault}{\updefault}$\delta\!-\!8$}}}}
\put(6398,-1071){\makebox(0,0)[lb]{\smash{{\SetFigFont{9}{10.8}{\rmdefault}{\mddefault}{\updefault}$4\!-\!\delta$}}}}
\put(5687,-801){\makebox(0,0)[lb]{\smash{{\SetFigFont{9}{10.8}{\rmdefault}{\mddefault}{\updefault}$4\al\!+\!8$}}}}
\put(7167,-801){\makebox(0,0)[lb]{\smash{{\SetFigFont{9}{10.8}{\rmdefault}{\mddefault}{\updefault}$4\al$}}}}
\put(6285,437){\makebox(0,0)[lb]{\smash{{\SetFigFont{9}{10.8}{\rmdefault}{\mddefault}{\updefault}$\delta\!-\!4\al\!-\!4$}}}}
\put(5827,-134){\makebox(0,0)[lb]{\smash{{\SetFigFont{9}{10.8}{\rmdefault}{\mddefault}{\updefault}$4\al\!+\!8$}}}}
\put(4042,-1618){\makebox(0,0)[lb]{\smash{{\SetFigFont{9}{10.8}{\rmdefault}{\mddefault}{\updefault}$\delta\!+\!4$}}}}
\put(4996,-1618){\makebox(0,0)[lb]{\smash{{\SetFigFont{9}{10.8}{\rmdefault}{\mddefault}{\updefault}$\delta\!-\!4$}}}}
\put(5306,-2288){\makebox(0,0)[lb]{\smash{{\SetFigFont{9}{10.8}{\rmdefault}{\mddefault}{\updefault}$\delta\!-\!4$}}}}
\put(5382,437){\makebox(0,0)[lb]{\smash{{\SetFigFont{10}{12.0}{\rmdefault}{\mddefault}{\updefault}${\m T}_-'$}}}}
\put(3875,-801){\makebox(0,0)[lb]{\smash{{\SetFigFont{9}{10.8}{\rmdefault}{\mddefault}{\updefault}$12$}}}}
\put(3299,-801){\makebox(0,0)[lb]{\smash{{\SetFigFont{9}{10.8}{\rmdefault}{\mddefault}{\updefault}$4$}}}}
\put(2999,-142){\makebox(0,0)[lb]{\smash{{\SetFigFont{9}{10.8}{\rmdefault}{\mddefault}{\updefault}$0$}}}}
\put(2188,-142){\makebox(0,0)[lb]{\smash{{\SetFigFont{9}{10.8}{\rmdefault}{\mddefault}{\updefault}$0$}}}}
\put(2536,437){\makebox(0,0)[lb]{\smash{{\SetFigFont{9}{10.8}{\rmdefault}{\mddefault}{\updefault}$4\delta$}}}}
\put(1801,437){\makebox(0,0)[lb]{\smash{{\SetFigFont{10}{12.0}{\rmdefault}{\mddefault}{\updefault}${\m T}_+'$}}}}
\put(2523,-251){\makebox(0,0)[lb]{\smash{{\SetFigFont{9}{10.8}{\rmdefault}{\mddefault}{\updefault}$0$}}}}
\put(3236,158){\makebox(0,0)[lb]{\smash{{\SetFigFont{10}{12.0}{\familydefault}{\mddefault}{\updefault}{\bf or}}}}}
\put(4606,-234){\makebox(0,0)[lb]{\smash{{\SetFigFont{9}{10.8}{\rmdefault}{\mddefault}{\updefault}$0$}}}}
\put(4475,-808){\makebox(0,0)[lb]{\smash{{\SetFigFont{9}{10.8}{\rmdefault}{\mddefault}{\updefault}$4\al\!-\!4$}}}}
\put(6462,-808){\makebox(0,0)[lb]{\smash{{\SetFigFont{9}{10.8}{\rmdefault}{\mddefault}{\updefault}$0$}}}}
\put(6418,-1756){\makebox(0,0)[lb]{\smash{{\SetFigFont{9}{10.8}{\rmdefault}{\mddefault}{\updefault}$4$}}}}
\put(6431,-2321){\makebox(0,0)[lb]{\smash{{\SetFigFont{9}{10.8}{\rmdefault}{\mddefault}{\updefault}$4\al$}}}}
\put(3836,-2309){\makebox(0,0)[lb]{\smash{{\SetFigFont{9}{10.8}{\rmdefault}{\mddefault}{\updefault}$\delta\!+\!4$}}}}
\put(4345,-2321){\makebox(0,0)[lb]{\smash{{\SetFigFont{9}{10.8}{\rmdefault}{\mddefault}{\updefault}$4\al\!+\!4\!\!-\!\delta$}}}}
\put(4496,-1712){\makebox(0,0)[lb]{\smash{{\SetFigFont{9}{10.8}{\rmdefault}{\mddefault}{\updefault}$8\!-\!\delta$}}}}
\put(6349,-216){\makebox(0,0)[lb]{\smash{{\SetFigFont{9}{10.8}{\rmdefault}{\mddefault}{\updefault}$4\!-\!4\al$}}}}
\put(1801,-821){\makebox(0,0)[lb]{\smash{{\SetFigFont{9}{10.8}{\rmdefault}{\mddefault}{\updefault}$4\al$}}}}
\put(2584,-821){\makebox(0,0)[lb]{\smash{{\SetFigFont{9}{10.8}{\rmdefault}{\mddefault}{\updefault}$4$}}}}
\end{picture}%

%% file: fig21y22.pstex_t
\begin{picture}(0,0)%
\includegraphics{fig21y22.pstex}%
\end{picture}%
\setlength{\unitlength}{3947sp}%
\begingroup\makeatletter\ifx\SetFigFont\undefined%
\gdef\SetFigFont#1#2#3#4#5{%
  \reset@font\fontsize{#1}{#2pt}%
  \fontfamily{#3}\fontseries{#4}\fontshape{#5}%
  \selectfont}%
\fi\endgroup%
\begin{picture}(5882,1241)(526,-846)
\put(4995,-810){\makebox(0,0)[lb]{\smash{{\SetFigFont{10}{12.0}{\rmdefault}{\mddefault}{\updefault}$8$}}}}
\put(5578,-810){\makebox(0,0)[lb]{\smash{{\SetFigFont{10}{12.0}{\rmdefault}{\mddefault}{\updefault}$0$}}}}
\put(6162,-810){\makebox(0,0)[lb]{\smash{{\SetFigFont{10}{12.0}{\rmdefault}{\mddefault}{\updefault}$0$}}}}
\put(5254,-227){\makebox(0,0)[lb]{\smash{{\SetFigFont{10}{12.0}{\rmdefault}{\mddefault}{\updefault}$8$}}}}
\put(5514,-318){\makebox(0,0)[lb]{\smash{{\SetFigFont{10}{12.0}{\rmdefault}{\mddefault}{\updefault}$4$}}}}
\put(5907,-232){\makebox(0,0)[lb]{\smash{{\SetFigFont{10}{12.0}{\rmdefault}{\mddefault}{\updefault}$0$}}}}
\put(5578,259){\makebox(0,0)[lb]{\smash{{\SetFigFont{10}{12.0}{\rmdefault}{\mddefault}{\updefault}$0$}}}}
\put(3259,189){\makebox(0,0)[lb]{\smash{{\SetFigFont{12}{14.4}{\rmdefault}{\mddefault}{\updefault}and}}}}
\put(526,-810){\makebox(0,0)[lb]{\smash{{\SetFigFont{10}{12.0}{\rmdefault}{\mddefault}{\updefault}$2$}}}}
\put(1109,-810){\makebox(0,0)[lb]{\smash{{\SetFigFont{10}{12.0}{\rmdefault}{\mddefault}{\updefault}$2$}}}}
\put(1693,-810){\makebox(0,0)[lb]{\smash{{\SetFigFont{10}{12.0}{\rmdefault}{\mddefault}{\updefault}$0$}}}}
\put(786,-227){\makebox(0,0)[lb]{\smash{{\SetFigFont{10}{12.0}{\rmdefault}{\mddefault}{\updefault}$0$}}}}
\put(1045,-318){\makebox(0,0)[lb]{\smash{{\SetFigFont{10}{12.0}{\rmdefault}{\mddefault}{\updefault}$0$}}}}
\put(1438,-236){\makebox(0,0)[lb]{\smash{{\SetFigFont{10}{12.0}{\rmdefault}{\mddefault}{\updefault}$0$}}}}
\put(2601,-810){\makebox(0,0)[lb]{\smash{{\SetFigFont{10}{12.0}{\rmdefault}{\mddefault}{\updefault}$4$}}}}
\put(2211,-227){\makebox(0,0)[lb]{\smash{{\SetFigFont{10}{12.0}{\rmdefault}{\mddefault}{\updefault}$0$}}}}
\put(2536,-292){\makebox(0,0)[lb]{\smash{{\SetFigFont{10}{12.0}{\rmdefault}{\mddefault}{\updefault}$0$}}}}
\put(2957,-260){\makebox(0,0)[lb]{\smash{{\SetFigFont{10}{12.0}{\rmdefault}{\mddefault}{\updefault}$0$}}}}
\put(3185,-810){\makebox(0,0)[lb]{\smash{{\SetFigFont{10}{12.0}{\rmdefault}{\mddefault}{\updefault}$4$}}}}
\put(1993,-810){\makebox(0,0)[lb]{\smash{{\SetFigFont{10}{12.0}{\rmdefault}{\mddefault}{\updefault}$0$}}}}
\put(2526,256){\makebox(0,0)[lb]{\smash{{\SetFigFont{10}{12.0}{\rmdefault}{\mddefault}{\updefault}$4\delta$}}}}
\put(993,256){\makebox(0,0)[lb]{\smash{{\SetFigFont{10}{12.0}{\rmdefault}{\mddefault}{\updefault}$16\delta$}}}}
\put(3570,-810){\makebox(0,0)[lb]{\smash{{\SetFigFont{10}{12.0}{\rmdefault}{\mddefault}{\updefault}$2$}}}}
\put(4154,-810){\makebox(0,0)[lb]{\smash{{\SetFigFont{10}{12.0}{\rmdefault}{\mddefault}{\updefault}$2$}}}}
\put(4738,-810){\makebox(0,0)[lb]{\smash{{\SetFigFont{10}{12.0}{\rmdefault}{\mddefault}{\updefault}$0$}}}}
\put(3797,-227){\makebox(0,0)[lb]{\smash{{\SetFigFont{10}{12.0}{\rmdefault}{\mddefault}{\updefault}$0$}}}}
\put(4090,-324){\makebox(0,0)[lb]{\smash{{\SetFigFont{10}{12.0}{\rmdefault}{\mddefault}{\updefault}$0$}}}}
\put(4511,-227){\makebox(0,0)[lb]{\smash{{\SetFigFont{10}{12.0}{\rmdefault}{\mddefault}{\updefault}$0$}}}}
\put(4059,256){\makebox(0,0)[lb]{\smash{{\SetFigFont{10}{12.0}{\rmdefault}{\mddefault}{\updefault}$64$}}}}
\put(4526,256){\makebox(0,0)[lb]{\smash{{\SetFigFont{12}{14.4}{\rmdefault}{\mddefault}{\updefault}${\m T}_+\ne {\m T'}_-$}}}}
\put(1459,256){\makebox(0,0)[lb]{\smash{{\SetFigFont{12}{14.4}{\rmdefault}{\mddefault}{\updefault}${\m T}_+\ne {\m T'}_+$}}}}
\end{picture}%

%% file: fig23.pstex_t
\begin{picture}(0,0)%
\includegraphics{fig23.pstex}%
\end{picture}%
\setlength{\unitlength}{3947sp}%
\begingroup\makeatletter\ifx\SetFigFont\undefined%
\gdef\SetFigFont#1#2#3#4#5{%
  \reset@font\fontsize{#1}{#2pt}%
  \fontfamily{#3}\fontseries{#4}\fontshape{#5}%
  \selectfont}%
\fi\endgroup%
\begin{picture}(5846,1241)(451,-846)
\put(1318,256){\makebox(0,0)[lb]{\smash{{\SetFigFont{12}{14.4}{\rmdefault}{\mddefault}{\updefault}${\m T}_-\ne {\m T'}_-$}}}}
\put(3138,162){\makebox(0,0)[lb]{\smash{{\SetFigFont{12}{14.4}{\rmdefault}{\mddefault}{\updefault}and}}}}
\put(3495,-810){\makebox(0,0)[lb]{\smash{{\SetFigFont{10}{12.0}{\rmdefault}{\mddefault}{\updefault}$2$}}}}
\put(4079,-810){\makebox(0,0)[lb]{\smash{{\SetFigFont{10}{12.0}{\rmdefault}{\mddefault}{\updefault}$4$}}}}
\put(4663,-810){\makebox(0,0)[lb]{\smash{{\SetFigFont{10}{12.0}{\rmdefault}{\mddefault}{\updefault}$2$}}}}
\put(3722,-227){\makebox(0,0)[lb]{\smash{{\SetFigFont{10}{12.0}{\rmdefault}{\mddefault}{\updefault}$2$}}}}
\put(4015,-324){\makebox(0,0)[lb]{\smash{{\SetFigFont{10}{12.0}{\rmdefault}{\mddefault}{\updefault}$0$}}}}
\put(4436,-227){\makebox(0,0)[lb]{\smash{{\SetFigFont{10}{12.0}{\rmdefault}{\mddefault}{\updefault}$2$}}}}
\put(4079,259){\makebox(0,0)[lb]{\smash{{\SetFigFont{10}{12.0}{\rmdefault}{\mddefault}{\updefault}$6$}}}}
\put(4451,256){\makebox(0,0)[lb]{\smash{{\SetFigFont{12}{14.4}{\rmdefault}{\mddefault}{\updefault}${\m T}_-\ne {\m T'}_-$}}}}
\put(967,-810){\makebox(0,0)[lb]{\smash{{\SetFigFont{10}{12.0}{\rmdefault}{\mddefault}{\updefault}$4$}}}}
\put(1551,-810){\makebox(0,0)[lb]{\smash{{\SetFigFont{10}{12.0}{\rmdefault}{\mddefault}{\updefault}$2$}}}}
\put(644,-227){\makebox(0,0)[lb]{\smash{{\SetFigFont{10}{12.0}{\rmdefault}{\mddefault}{\updefault}$2$}}}}
\put(903,-318){\makebox(0,0)[lb]{\smash{{\SetFigFont{10}{12.0}{\rmdefault}{\mddefault}{\updefault}$0$}}}}
\put(1296,-236){\makebox(0,0)[lb]{\smash{{\SetFigFont{10}{12.0}{\rmdefault}{\mddefault}{\updefault}$2$}}}}
\put(935,259){\makebox(0,0)[lb]{\smash{{\SetFigFont{10}{12.0}{\rmdefault}{\mddefault}{\updefault}$22$}}}}
\put(1811,-810){\makebox(0,0)[lb]{\smash{{\SetFigFont{10}{12.0}{\rmdefault}{\mddefault}{\updefault}$10$}}}}
\put(2459,-810){\makebox(0,0)[lb]{\smash{{\SetFigFont{10}{12.0}{\rmdefault}{\mddefault}{\updefault}$2$}}}}
\put(3043,-810){\makebox(0,0)[lb]{\smash{{\SetFigFont{10}{12.0}{\rmdefault}{\mddefault}{\updefault}$2$}}}}
\put(2070,-227){\makebox(0,0)[lb]{\smash{{\SetFigFont{10}{12.0}{\rmdefault}{\mddefault}{\updefault}$10$}}}}
\put(2394,-292){\makebox(0,0)[lb]{\smash{{\SetFigFont{10}{12.0}{\rmdefault}{\mddefault}{\updefault}2}}}}
\put(2815,-260){\makebox(0,0)[lb]{\smash{{\SetFigFont{10}{12.0}{\rmdefault}{\mddefault}{\updefault}$2$}}}}
\put(2459,259){\makebox(0,0)[lb]{\smash{{\SetFigFont{10}{12.0}{\rmdefault}{\mddefault}{\updefault}$0$}}}}
\put(4920,-810){\makebox(0,0)[lb]{\smash{{\SetFigFont{10}{12.0}{\rmdefault}{\mddefault}{\updefault}$4$}}}}
\put(5503,-810){\makebox(0,0)[lb]{\smash{{\SetFigFont{10}{12.0}{\rmdefault}{\mddefault}{\updefault}$4$}}}}
\put(5179,-227){\makebox(0,0)[lb]{\smash{{\SetFigFont{10}{12.0}{\rmdefault}{\mddefault}{\updefault}$4$}}}}
\put(5439,-318){\makebox(0,0)[lb]{\smash{{\SetFigFont{10}{12.0}{\rmdefault}{\mddefault}{\updefault}$4$}}}}
\put(5832,-232){\makebox(0,0)[lb]{\smash{{\SetFigFont{10}{12.0}{\rmdefault}{\mddefault}{\updefault}$0$}}}}
\put(5503,259){\makebox(0,0)[lb]{\smash{{\SetFigFont{10}{12.0}{\rmdefault}{\mddefault}{\updefault}$2$}}}}
\put(6051,-810){\makebox(0,0)[lb]{\smash{{\SetFigFont{10}{12.0}{\rmdefault}{\mddefault}{\updefault}$0$}}}}
\put(451,-810){\makebox(0,0)[lb]{\smash{{\SetFigFont{10}{12.0}{\rmdefault}{\mddefault}{\updefault}$2$}}}}
\end{picture}%

%% file: fig24.pstex_t
\begin{picture}(0,0)%
\includegraphics{fig24.pstex}%
\end{picture}%
\setlength{\unitlength}{3947sp}%
\begingroup\makeatletter\ifx\SetFigFont\undefined%
\gdef\SetFigFont#1#2#3#4#5{%
  \reset@font\fontsize{#1}{#2pt}%
  \fontfamily{#3}\fontseries{#4}\fontshape{#5}%
  \selectfont}%
\fi\endgroup%
\begin{picture}(5781,2027)(751,-2062)
\put(751,-1125){\makebox(0,0)[lb]{\smash{{\SetFigFont{14}{16.8}{\rmdefault}{\mddefault}{\updefault}$\x$}}}}
\put(1585,-181){\makebox(0,0)[lb]{\smash{{\SetFigFont{14}{16.8}{\rmdefault}{\mddefault}{\updefault}$\y$}}}}
\put(1731,-2016){\makebox(0,0)[lb]{\smash{{\SetFigFont{12}{14.4}{\rmdefault}{\mddefault}{\updefault}(norm $4C$)}}}}
\put(1520,-2016){\makebox(0,0)[lb]{\smash{{\SetFigFont{14}{16.8}{\rmdefault}{\mddefault}{\updefault}$\z$}}}}
\put(4740,-2016){\makebox(0,0)[lb]{\smash{{\SetFigFont{14}{16.8}{\rmdefault}{\mddefault}{\updefault}$\z$}}}}
\put(4951,-2016){\makebox(0,0)[lb]{\smash{{\SetFigFont{12}{14.4}{\rmdefault}{\mddefault}{\updefault}(norm $C$)}}}}
\put(3879,-1125){\makebox(0,0)[lb]{\smash{{\SetFigFont{14}{16.8}{\rmdefault}{\mddefault}{\updefault}$\x$}}}}
\put(4651,-211){\makebox(0,0)[lb]{\smash{{\SetFigFont{14}{16.8}{\rmdefault}{\mddefault}{\updefault}$\y$}}}}
\end{picture}%